\documentclass[11pt,twoside]{article}

\usepackage{amsbsy,amsfonts,amsmath,amssymb,eucal,mathrsfs}
\usepackage[all]{xy}
\usepackage[T1]{fontenc}
\usepackage[french]{babel}
\usepackage{tikz}
\usepackage{hyperref}

\def\itemn#1{\item[\hspace{0.6mm} {\rm (#1)}]}
\def\itemm#1{\item[\indent {\rm (#1)}]}
\def\no{\noindent}
\def\qed{\hfill $\square$ \vspace{5mm}}

\def\red{{\rm red}}

\def\cco{c.c.o.}
\def\cio{c.i.o.}
\def\ccf{c.c.f.}
\def\cif{c.i.f.}

\def\i{{\rm (i)}}

\def\into{\hookrightarrow}

\renewcommand{\tilde}{\widetilde}

\renewcommand{\bar}{\overline}

\newtheorem{counter}[subsubsection]{$\!\!$}
\newtheorem{subcounter}[subsection]{$\!\!$}
\newenvironment{defi}{\begin{counter} \rm {\bf D{\'e}finition.}}{\end{counter}}
\newenvironment{defis}{\begin{counter} \rm {\bf D{\'e}finitions.}}{\end{counter}}

\newenvironment{prop}{\begin{counter} {\bf Proposition.}}{\end{counter}}
\newenvironment{lemm}{\begin{counter} {\bf Lemme.}}{\end{counter}}
\newenvironment{coro}{\begin{counter} {\bf Corollaire.}}{\end{counter}}
\newenvironment{theo}{\begin{counter} {\bf Th{\'e}or{\`e}me.}}{\end{counter}}
\newenvironment{quot}{\begin{counter} {\bf Th{\'e}or{\`e}me }}{\end{counter}}

\newenvironment{rema}{\begin{counter} \rm {\bf Remarque.}}{\end{counter}}
\newenvironment{remas}{\begin{counter} \rm {\bf Remarques.}}{\end{counter}}

\newenvironment{exems}{\begin{counter} \rm {\bf Exemples.}}{\end{counter}}

\newenvironment{adoc}[1]{\begin{counter} \rm {\bf #1.}}{\end{counter}}
\newenvironment{demo}{{\flushleft \bf Preuve~:}}{\hfill $\square$ \vspace{5mm}}

\newenvironment{appdefi}{\begin{subcounter} \rm {\bf D{\'e}finition.}}
{\end{subcounter}}

\newenvironment{apptheo}{\begin{subcounter} {\bf Th{\'e}or{\`e}me.}}
{\end{subcounter}}

\newenvironment{applemm}{\begin{subcounter} {\bf Lemme.}}{\end{subcounter}}
\newenvironment{appadoc}[1]{\begin{subcounter} \rm {\bf #1.}}{\end{subcounter}}

\DeclareMathOperator{\Aut}{Aut}
\DeclareMathOperator{\Obj}{Obj}
\DeclareMathOperator{\Hom}{Hom}
\DeclareMathOperator{\Ext}{Ext}

\DeclareMathOperator{\Spec}{Spec}
\DeclareMathOperator{\Proj}{Proj}

\DeclareMathOperator{\SL}{SL}

\DeclareMathOperator{\disc}{disc}
\DeclareMathOperator{\Ass}{Ass}

\DeclareMathOperator{\Fitt}{Fitt}

\DeclareMathOperator{\Irr}{Irr}
\DeclareMathOperator{\St}{St}
\DeclareMathOperator{\Hilb}{Hilb}
\DeclareMathOperator{\Quot}{Quot}
\DeclareMathOperator{\Supp}{Supp}

   \def\cI{{\cal I}}
\def\cJ{{\cal J}}   
 \def\cO{{\cal O}}

\def\sC{{\mathscr C}} \def\sD{{\mathscr D}} \def\sF{{\mathscr F}}
\def\sG{{\mathscr G}} \def\sH{{\mathscr H}} \def\sI{{\mathscr I}}
 \def\sM{{\mathscr M}} \def\sP{{\mathscr P}}
\def\sR{{\mathscr R}} \def\sS{{\mathscr S}} \def\sU{{\mathscr U}}
\def\sV{{\mathscr V}} \def\sW{{\mathscr W}} \def\sX{{\mathscr X}}
\def\sY{{\mathscr Y}} \def\sZ{{\mathscr Z}}

\def\bP{{\mathbb P}} \def\bQ{{\mathbb Q}} \def\bZ{{\mathbb Z}}

\def\ff{{\mathfrak f}}  \def\fr{{\mathfrak r}}

\begin{document}

\begin{center}
{\bf \Large Composantes connexes et irr{\'e}ductibles en familles}

\bigskip
\bigskip

{\bf Matthieu Romagny}

\medskip

{\em 14 d{\'e}cembre 2009}
\end{center}

\bigskip

{\def\thefootnote{\relax}
\footnote{ \hspace{-6.8mm}
Mots cl{\'e}~: composantes connexes, composantes irr{\'e}ductibles,
crit{\`e}res d'Artin, sch{\'e}ma de Hilbert. 
Matthieu {\sc Romagny}, Institut de Math{\'e}matiques,
Th{\'e}orie des Nombres, Universit{\'e} Pierre et Marie Curie, Case~82,
4,~place Jussieu,
F-75252 Paris Cedex 05.

\noindent {\it email}~: romagny@math.jussieu.fr
}}

\begin{center}
\bf{Abstract}
\end{center}

For an algebraic stack $\sX$ flat and of finite presentation
over a scheme $S$, we introduce various notions of
{\em relative connected components} and {\em relative
irreducible components}. The main distinction between these
notions is whether we require the total space of a relative
component to be open or closed in $\sX$. We study the
representability of the associated functors of relative components,
and give an application to the moduli stack of curves of genus $g$
admitting an action of a fixed finite group $G$.

\section{Introduction}

\begin{appadoc}{Motivation}
La preuve de l'irr{\'e}ductibilit{\'e} de l'espace de modules des
courbes de genre $g$ par Deligne et Mumford en 1969 utilise
le fait que pour un morphisme $\sX\to S$ propre,
plat, de pr{\'e}sentation finie, {\`a} fibres g{\'e}om{\'e}triquement
normales, le nombre g{\'e}om{\'e}trique de composantes
irr{\'e}ductibles des fibres est constant. Dans cette situation
comme dans bien d'autres en g{\'e}om{\'e}trie alg{\'e}brique relative,
on voudrait en savoir un peu plus sur la variation des
composantes connexes et irr{\'e}ductibles dans les fibres d'une
famille. Ainsi, dans l'exemple ci-dessus, on peut penser qu'il
existe en fait un espace de modules pour les {\em composantes
irr{\'e}ductibles relatives} ({\`a} d{\'e}finir) de l'espace de
modules des courbes (lisses, ou stables) de genre $g$ qui est
repr{\'e}sentable par le sch{\'e}ma $\Spec(\bZ)$.
La contribution principale du pr\'esent article est de proposer
diff\'erentes notions de composantes connexes et irr{\'e}ductibles
en familles, de montrer que les foncteurs auxquelles elles donnent
naissance sont repr\'esentables, de comparer ces notions lorsque
c'est possible, et de donner des exemples et contre-exemples.
En application, nous \'etudions les composantes irr\'eductibles
du lieu $\sM_g(G)$ des courbes qui admettent une action d'un
groupe fini $G$ fix\'e.
\end{appadoc}

\begin{appadoc}{Composantes connexes}
Il est l\'eg\`erement plus simple de d\'evelopper la th\'eorie
des composantes connexes~; parlons donc d'abord de celles-ci.
Au minimum, une {\em composante connexe relative} pour
$\sX\to S$ doit \^etre un sous-champ $\sC\subset\sX$ qui est
plat et de pr\'esentation finie sur $S$, et dont chaque fibre
g\'eom\'etrique est une composante connexe de la fibre
correspondante de $\sX$. Il y a ensuite deux notions
naturelles possibles, selon que l'on demande {\`a} une telle
composante $\sC$ d'{\^e}tre ouverte ou ferm{\'e}e dans $\sX$.

Comme il s'agit d'une notion de nature essentiellement
topologique, on aimerait qu'une composante connexe relative
soit d{\'e}termin{\'e}e par son support. Ceci m{\`e}ne {\`a}
consid{\'e}rer en premier lieu des sous-champs {\em ouverts}
$\sC\subset \sX$.
Si $\sX$ est \`a fibres g\'eom\'etriquement r\'eduites, on
montre que le foncteur correspondant $\pi_0(\sX/S)$ est
repr\'esentable par un espace alg\'ebrique \'etale et
quasi-compact (th\'eor\`eme~\ref{theo:pi_0_representable}).
On obtient m\^eme une description tr\`es pr\'ecise de cet espace
comme quotient de $\sX$ par la relation d'\'equivalence d\'efinie
par l'appartenance \`a la m\^eme composante. Malheureusement,
dans le cas de fibres non g\'eom\'etriquement r\'eduites,
l'effectivit{\'e} des
composantes formelles est prise en d\'efaut et le foncteur des
composantes connexes n'est en g\'en\'eral pas repr\'esentable.

Alternativement, on peut aussi consid\'erer, en guise de
composantes relatives, des sous-champs {\em ferm\'es}. En
s'appuyant sur le sch\'ema de Hilbert, on montre que si
$\sX\to S$ est propre, le foncteur correspondant
$\pi_0(\sX/S)^\ff$ est repr\'esentable par un espace
alg\'ebrique formel, localement de pr\'esentation
finie et s\'epar\'e. Cet espace peut \^etre tr\`es (trop)
gros, mais le sous-foncteur $\pi_0(\sX/S)^\fr$
des composantes connexes ferm\'ees dont les fibres sont
g\'eom\'etriquement r\'eduites est repr\'esentable par un
sch\'ema formel quasi-fini et s\'epar\'e
(th\'eor\`eme~\ref{theo:pi_0_cas_propre}).

Lorsque $\sX\to S$ est \`a la fois \`a fibres
g\'eom\'etriquement r\'eduites et propre, on peut comparer
ces constructions. On montre que $\pi_0(\sX/S)$,
$\pi_0(\sX/S)^\ff$ et $\pi_0(\sX/S)^\fr$ sont alors
isomorphes, et ils sont aussi isomorphes \`a la
factorisation de Stein de $\sX$
(proposition~\ref{prop:propre_et_reduit}). On peut encore
les comparer lorsque $\sX\to S$ est seulement \`a fibres
g\'eom\'etriquement r\'eduites et pur~: on trouve que
$\pi_0(\sX/S)^\ff$ est un ouvert de $\pi_0(\sX/S)$.
\end{appadoc}

\begin{appadoc}{Composantes irr\'eductibles}
Les questions pos\'eees ci-dessus pour les composantes connexes
ont un analogue pour les composantes irr\'eductibles. Le
travail auquel cela m\`ene et les r\'eponses obtenues
semblent d'ailleurs un peu plus originaux.
On d\'efinit ainsi un foncteur
$\Irr(\sX/S)$ de composantes irr\'eductibles ouvertes, qui
est repr\'esentable par un espace alg\'ebrique \'etale et
quasi-compact si $\sX$ est \`a fibres g\'eom\'etriquement
r\'eduites (th\'eor\`eme~\ref{theo:pi_0_representable}).
On d\'efinit aussi un foncteur $\Irr(\sX/S)^\ff$ de
composantes irr\'eductibles ferm\'ees, mais on montre, en
regardant la famille donn\'ee par la conique plane universelle,
que m\^eme si $\sX$ est propre, ce foncteur n'est en
g\'en\'eral pas repr\'esentable par un un espace
alg\'ebrique formel (voir~\ref{adoc:contre_ex_coniques}).
Si $\sX$ est \`a fibres g\'eom\'etriquement r\'eduites,
ce foncteur est tout de m\^eme ouvert dans $\Irr(\sX/S)$.
\end{appadoc}


\begin{appadoc}{Application~: courbes avec action de $G$}
En guise d'application, nous d\'emontrons le
r\'esultat suivant. Soit $G$ un groupe fini, $\gamma$ son
cardinal, $g\ge 2$ un entier et $\sM_g(G)$ le sous-champ
du champ des courbes projectives lisses de genre $g$ form\'e
des courbes qui admettent une action fid{\`e}le de $G$.
Alors, sur le sch\'ema de base $S=\Spec(\bZ[1/30\gamma])$,
le foncteur des composantes irr\'eductibles ouvertes de
$\sM_g(G)$ est repr\'esentable par un sch\'ema fini \'etale
(corollaire~\ref{coro:Irr_de_MgG_fini_etale}). Si $G$ n'est
pas dans une liste explicite de $10$ groupes,
ce r\'esultat est m\^eme valable sur $\Spec(\bZ[1/2\gamma])$.
Notons que c'est l'int\'er\^et pour des objets tels que
$\sM_g(G)$ ou d'autres champs classifiants qui justifie
l'effort fait pour \'etablir les r\'esultats dans le cadre
des champs alg\'ebriques.
\end{appadoc}

\begin{appadoc}{Remarques}
Dans le cas o{\`u} $\sX$ est un $S$-sch{\'e}ma lisse et
quasi-compact, on trouve une br{\`e}ve {\'e}tude de $\pi_0(\sX/S)$
dans \cite{LMB},~(6.8). Son introduction dans {\em loc. cit.}
est motiv{\'e}e par l'int{\'e}r{\^e}t pour la notion
d'{\em {\'e}quiconnexit{\'e}}, notion {\'e}troitement li{\'e}e {\`a}
son tour aux propri{\'e}t{\'e}s de s{\'e}paration de $\pi_0(\sX/S)$.

De nombreux \'enonc\'es classiques de g\'eom\'etrie alg\'ebrique
connus pour les sch\'emas s'\'etendent, plus ou moins facilement
d'ailleurs, au cadre des champs alg\'ebriques. Il est d'usage
d'admettre purement et simplement les \'enonc\'es dont la preuve
est essentiellement la m\^eme que pour les sch\'emas, et j'avoue
c\'eder parfois \`a ce travers. Cependant, j'ai pr\'ef\'er\'e
v\'erifier soigneusement que les \'enonc\'es de \cite{EGA} sur
la constructiblit\'e de certaines parties et de certaines
propri\'et\'es s'adaptent bien aux champs, car ces \'enonc\'es
sont utilis\'es abondamment dans ce texte. J'ai donc inclus en
fin d'article un premier appendice sur ce sujet, et un second
appendice qui utilise le premier pour \'etendre aux champs
certains r\'esultats de \cite{Ro2} sur les sch\'emas purs,
utilis\'es dans le pr\'esent article.
\end{appadoc}

\begin{appadoc}{Notation}
Dans tout le texte, nous
utilisons les notation $(R,K,k,\pi)$ pour d\'esigner un anneau
de valuation discr\`ete $R$ de corps de fractions $K$, de
corps r\'esiduel $k$, et une uniformisante $\pi$.
\end{appadoc}

\begin{appadoc}{Remerciements}
En premier lieu, je souhaite remercier chaleureusement Pierre
Lochak qui m'a pos\'e une question \`a l'origine de ce travail,
en a suivi la progression et a \'ecout\'e mes r\'eflexions
sur la question depuis quelques mois. Je veux aussi remercier
Sylvain Brochard, qui m'a fourni l'une des deux preuves donn\'ees
ici d'un \'enonc\'e sur les composantes le long d'une section
dans le cadre des champs. J'exprime enfin
ma gratitude \`a Angelo Vistoli qui m'a sugg\'er\'e l'exemple
de la conique universelle, et \`a Sylvain Maugeais pour quelques
\'echanges de courriels sur le th\`eme de l'article.
\end{appadoc}

\tableofcontents

\section{Composantes ouvertes}

Une fois d\'efinies les composantes ouvertes, la preuve
de la repr\'esentabilit\'e des foncteurs associ\'es n\'ecessite
certains r\'esultats interm\'ediaires sur les composantes
connexes et irr\'eductibles le long d'une section. Ces
r\'esultats sont discut\'es dans~\ref{ss_section:cc_ci__enonces},
\ref{ss_section:cc_preuves} et~\ref{ss_section:ci_preuve}.
La repr\'esentabilit\'e et ses corollaires sont \'etablis
dans~\ref{ss_section:representabilite_de_pi_0}.

\subsection{D{\'e}finitions et remarques pr{\'e}liminaires}

Une composante connexe relative sera d\'efinie comme un
sous-champ $\sC\subset \sX$ plat sur la base et dont les
fibres g{\'e}om{\'e}triques sont des composantes connexes
des fibres correspondantes de $\sX$. Comme il s'agit d'une
notion essentiellement topologique, on aimerait qu'une
composante connexe relative soit d{\'e}termin{\'e}e par son
support. Ceci m{\`e}ne {\`a} consid\'erer des immersions
$\sC\subset \sX$ ouvertes.

Pour un espace topologique poss{\'e}dant un nombre
fini de composantes irr{\'e}ductibles, appelons {\em composante
irr{\'e}ductible ouverte} l'int{\'e}rieur d'une composante irr{\'e}ductible,
ou de mani{\`e}re {\'e}quivalente, le compl{\'e}mentaire
de toutes les composantes irr{\'e}ductibles sauf une.
Pour les m{\^e}mes raisons que ci-dessus, il est naturel de souhaiter
qu'une composante irr{\'e}ductible relative soit un sous-champ ouvert
$\sI\subset\sX$. Dans les fibres g{\'e}om{\'e}triques de $\sX$, on est
amen{\'e} {\`a} porter l'attention sur les composantes irr{\'e}ductibles
{\em ouvertes}. On arrive aux d{\'e}finitions suivantes.

\begin{defis} \label{defi:cco_et_cio}
Soit $\sX$ un champ alg{\'e}brique de pr{\'e}sentation finie sur un
sch{\'e}ma $S$.
\begin{trivlist}
\itemn{1} Une {\em composante connexe ouverte
(en abr{\'e}g{\'e} c.c.o.) de $\sX$ sur $S$}
est un sous-champ ouvert $\sC\subset\sX$ fid{\`e}lement plat
et de pr{\'e}sentation finie sur~$S$, tel que pour tout point
g{\'e}om{\'e}trique $\bar s:\Spec(\Omega)\to S$, la fibre
$\sC_{\bar s}$ est une composante connexe de $\sX_{\bar s}$.
On note $\pi_0(\sX/S)$ le foncteur qui {\`a} un $S$-sch{\'e}ma
$T$ associe l'ensemble des \cco\ de $\sX_T$ sur $T$.
\itemn{2} Une {\em composante irr{\'e}ductible ouverte
(en abr{\'e}g{\'e} c.i.o.) de $\sX$ sur $S$} est un sous-champ
ouvert $\sI\subset\sX$ fid{\`e}lement plat et de pr{\'e}sentation
finie sur~$S$, tel que pour tout point g{\'e}om{\'e}trique
$\bar s:\Spec(\Omega)\to S$, la fibre $\sI_{\bar s}$ est une
composante irr{\'e}ductible ouverte de $\sX_{\bar s}$. On note
$\Irr(\sX/S)$ le foncteur qui {\`a} un $S$-sch{\'e}ma $T$
associe l'ensemble des \cio\ de $\sX_T$ sur $T$.
\end{trivlist}
\end{defis}


La formation de ces foncteurs commute aux changements
de base $S'\to S$.

\begin{lemm} \label{lemme:faisceau_etale}
Soit $\sX$ un $S$-champ alg{\'e}brique plat et de pr{\'e}sentation
finie et soit $F$ l'un des deux foncteurs $\pi_0(\sX/S)$ ou
$\Irr(\sX/S)$. Alors, $F$ est un faisceau pour la topologie
{\'e}tale, {\`a} diagonale ouverte quasi-compacte. De plus $F$ est
{\'e}tale et quasi-compact sur $S$.
\end{lemm}

On rappelle que par d\'efinition, un faisceau est {\'e}tale
s'il est formellement
{\'e}tale et localement de pr{\'e}sentation finie.

\begin{demo}
La propri{\'e}t{\'e} de faisceau r{\'e}sulte
de faits classiques de th{\'e}orie de la descente. Pour voir que la
diagonale de $F$ est repr{\'e}sentable par une immersion ouverte
quasi-compacte, il suffit d'observer que si $T$ est un
$S$-sch{\'e}ma et $\sC,\sC'\in F(T)$, le lieu des points $t\in T$
tels que $\sC_t=\sC'_t$ est l'ouvert image de $\sC\cap \sC'$
dans~$T$. De plus, utilisant~\cite{EGA}~IV.8.6.3,
on voit que $F$ est localement de pr{\'e}sentation finie. Enfin,
si $T_0\to T$ une immersion ferm{\'e}e de sch{\'e}mas qui est un
hom{\'e}omorphisme, alors $\sX\times_S T_0\to \sX\times_S T$
est un hom{\'e}omorphisme, de sorte que $F(T)\to F(T_0)$ est bijectif.
Ceci montre que $F$ est formellement {\'e}tale sur $S$. Il ne reste
qu'{\`a} montrer que $F$ est quasi-compact sur $S$, ce qui r{\'e}sulte du
lemme~\ref{lemme:rep_sur_un_ouvert} ci-dessous.
\end{demo}

\begin{lemm} \label{lemme:rep_sur_un_ouvert}
Sous les m{\^e}mes hypoth{\`e}ses que dans~\ref{lemme:faisceau_etale},
il existe un ouvert $U$ contenant les points maximaux de $S$ tel que
la restriction de $F$ {\`a} $U$ soit repr{\'e}sentable par un $U$-espace
alg{\'e}brique quasi-compact.
\end{lemm}

On d{\'e}duit facilement de cet {\'e}nonc{\'e} qu'il existe une stratification
$S^*=\{S_i\}$ de $S$ telle que $F\times_S S^*$ est repr{\'e}sentable
par un $S^*$-espace alg{\'e}brique quasi-compact. En effet,
{\'e}tant donn{\'e} que $\sX$ est de pr{\'e}sentation finie sur $S$ et que la
formation de $F$ est compatible au changement de base, on peut se
ramener au cas o{\`u} $S$ est noeth{\'e}rien, auquel cas l'assertion provient
de~\ref{lemme:rep_sur_un_ouvert} par r{\'e}currence noeth{\'e}rienne. Puisque
$S^*\to S$ est bijectif, ceci prouve la quasi-compacit{\'e} annonc{\'e}e
dans~\ref{lemme:faisceau_etale}.

\begin{demo}
Le raisonnement {\'e}tant le m{\^e}me pour les deux foncteurs
consid{\'e}r{\'e}s, disons que $F=\pi_0(\sX/S)$ pour fixer les id{\'e}es.
Soit $\eta$ un point maximal de $S$. D'apr{\`e}s l'{\'e}nonc{\'e} sur
la diagonale dans le lemme~\ref{lemme:faisceau_etale}, il suffit de
montrer qu'il existe un voisinage ouvert de $\eta$ au-dessus duquel $F$
poss\`ede une pr\'esentation $P\to F$, i.e. un morphisme surjectif et
fppf (fid\`element plat de pr\'esentation finie) depuis un sch\'ema
quasi-compact $P$.
Comme $F$ est localement de pr\'esentation finie, on se ram\`ene au cas
o\`u $S$ est affine et de type fini sur $\bZ$. Comme $F$ est \'etale,
il revient au m\^eme de trouver une pr\'esentation au-dessus de $S$ ou
au-dessus de $S_\red$ donc on peut supposer $S$ r\'eduit. Comme $\eta$
est maximal, son anneau local est alors int\`egre donc quitte \`a
r\'etr\'ecir $S$, on peut le supposer int\`egre.

Notons $K$ une extension finie de $k(\eta)$ telle que les
sous-sch{\'e}mas ferm{\'e}s r{\'e}duits des composantes connexes
de $\sX_\eta\otimes K$ soient g{\'e}om{\'e}triquement connexes
et g\'eom\'etriquement r{\'e}duits~; on peut prendre pour $K$ un
compositum des corps de d{\'e}finition des composantes connexes
r\'eduites de $\sX_\eta\otimes \bar{k(\eta)}$.
Soit $S'$ un $S$-sch{\'e}ma affine int{\`e}gre et de type fini sur
$\bZ$, de corps de fonctions~$K$. Par platitude g{\'e}n{\'e}rique,
quitte \`a restreindre $S$ le morphisme $S'\to S$ est fppf. Alors
son image contient un ouvert, donc quitte {\`a} restreindre
encore~$S$ on peut supposer $S'\to S$ surjectif. Il suffit de
montrer l'assertion apr{\`e}s le changement de base $S'\to S$,
on peut donc remplacer $S$ par $S'$ et supposer que les
composantes connexes $\sC_{1,\eta}\dots,\sC_{r,\eta}$ de $\sX_\eta$
sont g{\'e}om{\'e}triquement connexes. Soient $\sC_1,\dots,\sC_r$
des sous-sch{\'e}mas ferm{\'e}s de $\sX$ qui induisent les
$\sC_{i,\eta}$. D'apr{\`e}s~\ref{prop:EGA_95}~(iv) et
\ref{theo:EGA_977}~(ii), quitte {\`a} r\'etr\'ecir encore $S$
on peut supposer que les $\sC_i$ sont les composantes connexes
relatives de $\sX$. Elles fournissent donc un morphisme surjectif
$S\times\{1,\dots,r\}\to F$, o{\`u} {\`a} la source la
notation d{\'e}signe la somme disjointe de $r$ copies de $S$.
Comme la source et le but sont \'etales, ce morphisme est \'etale,
et fournit donc une pr\'esentation de $\pi_0(\sX/S)$, qui est
donc un espace alg\'ebrique quasi-compact au-dessus
d'un voisinage de $\eta$.
\end{demo}

\begin{prop} \label{prop:pi_0_sur_un_corps}
Soit $\sX$ un champ alg{\'e}brique de type fini sur un corps $k$.
\begin{trivlist}
\itemn{1} Soit $A$ la plus grande sous-$k$-alg{\`e}bre s{\'e}parable
de $H^0(\sX,\cO_\sX)$. Alors $\pi_0(\sX/k)\simeq \Spec(A)$.
\itemn{2} Soit $B$ la cl{\^o}ture s{\'e}parable de $k$ dans l'anneau
total des fractions de $\sX$. Alors $\Irr(\sX/k)\simeq \Spec(B)$.
\end{trivlist}
\end{prop}

\begin{demo}
(1) D'apr{\`e}s le lemme pr{\'e}c{\'e}dent $\pi_0(\sX/k)$ est
repr{\'e}sentable par un $k$-espace alg{\'e}brique {\'e}tale et quasi-compact,
donc par un sch{\'e}ma affine. Par ailleurs, on dispose d'un morphisme
$f:\sX\to\Spec(A)$ {\`a} fibres g{\'e}om{\'e}triquement connexes
(voir \cite{EGA}~IV.4.5.15 qui est {\'e}nonc{\'e} pour les sch{\'e}mas mais
dont la preuve fonctionne {\`a} l'identique pour les champs).
Consid{\'e}rons $\sX$, vu comme $A$-sch{\'e}ma via $f$,
et le morphisme $\sX\to \sX\otimes_k A$ qui est une immersion ouverte
puisque $A$ est {\'e}tale sur $k$. Ce morphisme fait de $\sX$ une
\cco\ de $\sX\otimes_k A$ sur $\Spec(A)$, ce qui d{\'e}finit un morphisme
$g:\Spec(A)\to\pi_0(\sX/k)$. Comme~$g$ est un isomorphisme
apr{\`e}s passage {\`a} une cl{\^o}ture alg{\'e}brique de $k$, donc c'est un
isomorphisme.

\no (2) Notons $\sY$ la normalisation de $\sX_\red$. Nous anticipons
un peu sur des r{\'e}sultats de fonctorialit{\'e} qui seront
{\'e}tablis plus tard (voir
sous-section~\ref{subsection:fonctorialite}). On a
$\Irr(\sX/k)\simeq\Irr(\sX_\red/k)$, et comme le morphisme de
normalisation $\sY\to\sX_\red$ est birationnel, on a
$\Irr(\sX_\red/k)\simeq\Irr(\sY/k)$ (voir
corollaire~\ref{coro:covariance_en_X}). Par normalit{\'e}
on a $\Irr(\sY/k)=\pi_0(\sY/k)$ qui est repr{\'e}sentable par le
sch{\'e}ma fini, spectre de la plus grande sous-$k$-alg{\`e}bre
s{\'e}parable de $H^0(\sY,\cO_\sY)$
(proposition~\ref{prop:pi_0_sur_un_corps}), qui est aussi la plus
grande sous-$k$-alg{\`e}bre s{\'e}parable du corps de fonctions
$k(\sY)=k(\sX_\red)$. Le r{\'e}sultat en d{\'e}coule.
\end{demo}

\subsection{Composantes connexes et irr{\'e}ductibles
le long d'une section~: {\'e}nonc{\'e}s}
\label{ss_section:cc_ci__enonces}

Rappelons qu'un {\em point} d'un champ $\sX$ est une classe
d'{\'e}quivalence de points $x_K:\Spec(K)\to\sX$ {\`a} valeurs dans
un corps $K$, pour la relation qui identifie $x_K$ et $x_L$ si
et seulement s'il existe une extension $M$ de $K$ et $L$ telle
que les points $x_K\in\sX(K)$ et $x_L\in\sX(L)$ sont isomorphes
dans $\sX(M)$. On note $x=[x_K]$ le point ainsi d{\'e}fini et
$|\sX|$ l'espace topologique des points de $\sX$. Il y a une
bijection entre sous-ensembles ouverts de $|\sX|$ et
sous-champs ouverts de~$\sX$, et nous confondrons les deux.

\begin{prop} \label{comp_connexe_d_une_section}
Soit $\sX$ un $S$-champ alg{\'e}brique plat, de pr{\'e}sentation
finie, {\`a} fibres g{\'e}om{\'e}\-triquement r{\'e}duites et soit
$g:S\to\sX$ une section. Pour tout $s\in S$, on note $\sC_s$
la composante connexe de $g(s)$ dans $\sX_s$. Alors, la
r{\'e}union des $|\sC_s|$ est un ouvert $C\subset |\sX|$,
correspondant {\`a} un sous-champ ouvert $\sC\subset\sX$,
de pr{\'e}sentation finie sur $S$, et dont la formation commute
au changement de base.
\end{prop}

Nous appellerons $\sC$ la {\em \cco\ de $\sX$ le long de la
section $g$}. S'il est utile de pr{\'e}ciser la section, nous
noterons $\sC_s(g)$ au lieu de $\sC_s$ et $\sC(g)$ au lieu
de $\sC$.

\begin{defi}
Soit $X$ un espace topologique. On appelle {\em lieu
unicomposante de $X$} l'ensemble des points de $X$ qui
n'appartiennent qu'{\`a} une composante irr{\'e}ductible de $X$.
\end{defi}

Soit $\sX$ un $S$-champ alg{\'e}brique et $x$ un point de l'espace
topologique $|\sX|$. Notons $s$ l'image de $x$ dans $S$ et
d{\'e}signons par un indice $(\cdot)_{\bar s}$ le changement
de base au spectre d'une cl{\^o}ture alg{\'e}brique du corps
r{\'e}siduel de $s$. Il est facile de voir que la propri{\'e}t{\'e}
pour $x_{K,\bar s}:\Spec(K)_{\bar s}\to \sX_{\bar s}$ de se
factoriser par le lieu unicomposante de la fibre
g{\'e}om{\'e}trique $\sX_{\bar s}$ est ind{\'e}pendante du
repr{\'e}sentant $x_K:\Spec(K)\to\sX$ choisi pour $x$. Ceci
justifie la d{\'e}finition suivante.

\begin{defi}
Soit $\sX$ un $S$-champ alg{\'e}brique. On appelle
{\em lieu unicomposante (relatif) de $\sX$ sur~$S$} l'ensemble
des points $x\in|\sX|$, d'image $s$ dans $S$, tels que
$x_{K,\bar s}:\Spec(K)_{\bar s}\to \sX_{\bar s}$ se factorise
par le lieu unicomposante de $\sX_{\bar s}$, pour un (ou, de
mani{\`e}re {\'e}quivalente, n'importe quel)
repr{\'e}sentant $x_K:\Spec(K)\to\sX$.
\end{defi}

Noter que si $\sX$ est un champ sur un corps $k$, l'inclusion
du lieu unicomposante (relatif) de $\sX$ sur~$k$ dans le lieu
unicomposante (absolu) de $|\sX|$ est en g{\'e}n{\'e}ral stricte.
C'est le cas par exemple pour la conique sur $\bQ$ d'{\'e}quation
$x^2-2y^2=0$.

\begin{prop} \label{comp_irred_d_une_section}
Soit $\sX$ un $S$-champ alg{\'e}brique de pr{\'e}sentation
finie, {\`a} fibres g{\'e}om{\'e}triquement r{\'e}duites.
\begin{trivlist}
\itemn{i} Le lieu unicomposante de $\sX$ sur $S$ est un
ouvert $U\subset |\sX|$, correspondant {\`a} un sous-champ ouvert
$S$-dense $\sU\subset\sX$, de pr{\'e}sentation finie sur~$S$,
dont la formation commute au changement de base.
\itemn{ii} Soit $g:S\to \sX$ une section {\`a} valeurs dans $\sU$.
Pour tout $s\in S$, on note $\sI_s$ la composante
irr{\'e}ductible ouverte de $g(s)$ dans $\sX_s$. Alors, la r{\'e}union
des $|\sI_s|$ est un ouvert $I\subset |\sX|$, correspondant
{\`a} un sous-champ ouvert $\sI\subset \sX$, de pr{\'e}sentation
finie sur $S$, dont la formation commute au changement de base.
\end{trivlist}
\end{prop}

Nous appellerons $\sI$ la {\em \cio\ de $\sX$ le long
de la section $g$}. S'il est utile de pr{\'e}ciser la section,
nous noterons $\sI_s(g)$ au lieu de $\sI_s$ et $\sI(g)$
au lieu de $\sI$.

\begin{rema} \label{rema:necessite_de_la_platitude}
Le lemme~\ref{comp_connexe_d_une_section} est faux en g{\'e}n{\'e}ral
sans hypoth{\`e}se de platitude. Pour un contre-exemple, prenons
un anneau de valuation discr\`ete $(R,K,k,\pi)$ et le sch{\'e}ma
affine d'anneau $R[x,y]/(xy(y-1),\pi y(y-1))$ muni de n'importe
quelle $R$-section adh{\'e}rence d'un point $K$-rationnel de la fibre
g{\'e}n{\'e}rique. Il s'agit simplement du sch{\'e}ma obtenu
{\`a} partir du $R$-sch{\'e}ma plat r{\'e}union de trois droites
$xy(y-1)=0$, en enlevant de la fibre g{\'e}n{\'e}rique la composante
irr{\'e}ductible ouverte centrale d'{\'e}quation $x=0$.

Il est remarquable qu'au contraire, la platitude n'est pas
n{\'e}cessaire pour le lemme~\ref{comp_irred_d_une_section}.
\end{rema}

\subsection{Composantes connexes le long d'une section~: preuve}
\label{ss_section:cc_preuves}

Dans cette section, nous d{\'e}montrons la
proposition~\ref{comp_connexe_d_une_section}.
Nous d{\'e}montrons d'abord que la formation de la composante
connexe le long de la section $g$ commute au changement de
base. Ceci est cons{\'e}quence du lemme suivant~:

\begin{lemm} \label{lemm:geom_connexe}
Soient $k$ un corps et $f:\sY\to\sX$ un morphisme de
$k$-champs alg{\'e}briques. On suppose que $\sX$ est
connexe et $\sY$ est g{\'e}om{\'e}triquement connexe.
Alors, $\sX$ est g{\'e}om{\'e}triquement connexe.
\end{lemm}

\begin{demo}
La preuve de \cite{EGA}~IV.4.5.13 fonctionne {\`a} l'identique
pour les champs alg{\'e}briques.
\end{demo}

Pour tout morphisme de changement de base $S'\to S$, notons
$\sX'=\sX\times_S S'$, $g':S'\to\sX'$ la section d{\'e}duite
de la section $g$, et $C'$ la r{\'e}union des $|\sC_{s'}(g')|$.
Le fait que la formation de $\sC$ commute
au changement de base sera cons{\'e}quence du fait que $C'$
est la pr{\'e}image de $C$ par l'application continue
$|\sX'|\to |\sX|$. Ceci peut {\^e}tre v{\'e}rifi{\'e} fibre
{\`a} fibre. Or, si $s'\in S'$ a pour image $s\in S$,
la composante connexe $\sC_s$ est g{\'e}om{\'e}triquement
connexe d'apr{\`e}s le lemme~\ref{lemm:geom_connexe} appliqu{\'e}
{\`a} la section $g$. Il s'ensuit que
$\sC'_{s'}=\sC_s\otimes_{k(s)}k(s')$, d'o{\`u} notre assertion.

Il reste {\`a} d{\'e}montrer que $C$ est un ouvert.
Ce r{\'e}sultat est d{\'e}montr{\'e} pour les
sch{\'e}mas dans \cite{EGA}~IV.15.6.5. Nous donnerons deux
preuves de son extension au cas des champs alg{\'e}briques.
La premi{\`e}re m'a {\'e}t{\'e} communiqu{\'e}e par Sylvain
Brochard et proc{\`e}de par r{\'e}duction au cas des sch{\'e}mas.
La deuxi{\`e}me suit directement la preuve de \cite{EGA} pour
les sch{\'e}mas~; sa raison d'{\^e}tre est que, comme nous
l'expliquons dans~\ref{adoc:remarques}, elle s'adaptera
mieux que la premi{\`e}re au cas des composantes
irr{\'e}ductibles le long d'une section, et elle servira alors
de mod{\`e}le.

\begin{adoc}{Premi{\`e}re preuve} \label{preuve_de_Brochard}
On proc{\`e}de par r{\'e}duction au cas des sch{\'e}mas,
avec les m{\'e}thodes du paragraphe 4.2 et notamment du lemme
4.2.8 de~\cite{Bro}. Dans {\em loc. cit.}, la notation pour
$\sC$ est $X^0$ et pour faciliter la transcription, nous
adopterons cette notation dans le reste
de~\ref{preuve_de_Brochard}. On note $X=|\sX|$ et on adopte
les notations du diagramme suivant~:
$$
\xymatrix{
S'_1 \ar[r]^{\pi'} & \ar @{} [dr] |{\square} S_1 \ar[r]^{e_1} \ar[d]
& X_1 \ar[d]^\pi \ar@/^1.5pc/[dd]^{f_1=f\circ \pi} \\
& S \ar[r]^e \ar@{=}[rd] & \sX \ar[d]^f \\
& & S}
$$
dans lequel $\pi:X_1\to\sX$ est une pr{\'e}sentation de $\sX$
par un sch{\'e}ma et $\pi':S'_1\to S_1$ est une pr{\'e}sentation
de l'espace alg{\'e}brique $S_1:=S\times_\sX X_1$ par un
sch{\'e}ma. On consid{\`e}re le diagramme de sch{\'e}mas~:
$$
\xymatrix{S'_1 \ar[r]^{e_1\circ \pi'} & X_1 \ar[r]^{f_1} & S .}
$$
On note $W_0:=X^0_1(e_1\circ \pi')$ ({\em cf.} lemme 4.2.7
de \cite{Bro}) le sous-ensemble de $X_1$ dont la fibre
au-dessus d'un point $s\in S$ est la r{\'e}union des composantes
connexes de $(X_1)_s$ qui rencontrent $e_1\circ\pi'(S'_1)$.
Comme l'image continue d'un connexe est connexe, on voit que
$\pi(W_0)$ est inclus dans $X^0$. Vu que $\pi$ est lisse,
on en d{\'e}duit que $f_1$ est universellement ouvert en tout
point de $W_0$. Par \cite{Bro}, 4.2.7 (ii) a) il en r{\'e}sulte
que $W_0$ est un ouvert de $X_1$.

On consid{\`e}re maintenant $V_1=\pi^{-1}(\pi(W_0))$. C'est
un ouvert de $X_1$ qui contient $W_0$ (c'est le satur{\'e}
de $W_0$ pour la relation d'{\'e}quivalence d{\'e}finie par
$\pi$). D'apr{\`e}s~\cite{Bro}, 4.2.7 (ii) b) appliqu{\'e}
au diagramme
$$
\xymatrix{V_1 \ar@{^(->}[r]^{e_1\circ\pi'}
& X_1 \ar[r]^{f_1} & S}
$$
on voit que $W_1:=X^0_1(V_1\into X_1)$ est un ouvert de $X_1$.
On poursuit le processus en posant successivement
$$
V_2=\pi^{-1}(\pi(W_1)) \ ; \
W_2=X^0_1(V_2\into X_1) \ ; \
\dots \ ; \
V_n=\pi^{-1}(\pi(W_{n-1})) \ ; \
W_n=X^0_1(V_n\into X_1) \dots
$$
{\`A} chaque fois, on a $\pi(W_i)\subset X^0$ si bien que $f_1$
est universellement ouvert en tout point de $W_i$ et le lemme
4.2.7 (ii) b) de \cite{Bro} permet d'affirmer que $W_i$ est
ouvert. On obtient ainsi une suite croissante d'ouverts de
$X_1$~:
$$
W_0\subset V_1\subset W_1\subset V_2\subset W_2\subset\dots
$$
On note $V=\cup_i V_i=\cup_i W_i$. Pour conclure, il suffit
de montrer que $\pi(V)=X^0$. On peut pour cela raisonner
fibre par fibre et supposer que $S$ est le spectre d'un
corps. Quitte {\`a} remplacer $X$ par $X^0$, on peut aussi
supposer que $X$ est connexe (c'est-{\`a}-dire que $X^0=X$).
Vu la construction des $V_i$ (resp. des $W_i$) il est clair
que $V$ est une r{\'e}union de fibres de $\pi$ (resp. de
composantes connexes de $X_1$).

Pour conclure, il suffit de montrer que $\pi(V)$ contient toutes
les composantes irr{\'e}ductibles qu'il rencontre, car il sera
alors ferm{\'e}, ouvert et non vide dans $X$ connexe. Soient donc
$x\in\pi(V)$ et $Z$ une composante irr{\'e}ductible de $X$ qui
contient $x$. Soit $z\in Z$, il faut montrer que $z\in\pi(V)$.
Soit $U$ un ouvert connexe de $X_1$ dont l'image contient $z$.
Comme $\pi(V)\cap Z$ et $\pi(U)\cap Z$ sont deux ouverts non
vides de $Z$, leur intersection est non vide. Comme par
ailleurs $V$ est une r{\'e}union de fibres de $\pi$, on en
d{\'e}duit que $U\cap V$ est non vide. Enfin, vu que $U$ est
connexe et que $V$ est une r{\'e}union de composantes connexes
de $X_1$ on voit que $U$ est inclus dans $V$, si bien que
$z\in \pi(V)$. Ceci conclut la preuve
de~\ref{comp_connexe_d_une_section}.
\end{adoc}

\begin{adoc}{Deuxi{\`e}me preuve} \label{preuve_classique}
Revenons aux notations $\sC(g)$, $C(g)$ pour d\'esigner les
composantes connexes le long d'une section.
Nous suivons la preuve de~\cite{EGA}~IV.15.6.5, en soulignant
les modifications n{\'e}cessaires pour passer aux champs
alg{\'e}briques. Comme la formation de $C$ commute au changement
de base, en particulier on peut localiser et supposer $S=\Spec(A)$
affine. Par les arguments habituels, comme $\sX$ est de
pr{\'e}sentation finie sur $S$, on se ram{\`e}ne au cas o{\`u}
$A$ est noeth{\'e}rien.

On montre ensuite que $C$ est localement constructible. Les
arguments de \cite{EGA}~IV.9.7.12 restent
valables, avec quelques modifications mineures pour adapter aux
champs alg{\'e}briques les lemmes n{\'e}cessaires de \cite{EGA}~IV.9
sur la constructiblit{\'e}. Ces modifications sont indiqu{\'e}es dans
l'annexe~\ref{Annexe_constructiblite}, voir notamment
\ref{prop:EGA_95}~(i), \ref{prop:EGA_95}~(iv) et \ref{theo:EGA_977}~(ii).

Il reste {\`a} montrer que $C$ est stable par g{\'e}n{\'e}risation.
Pour cela, on est ramen{\'e} au cas o{\`u} $S$ est le spectre d'un
anneau de valuation discr{\`e}te $R=(R,K,k,\pi)$, et on peut
supposer $R$ complet et $k$ alg{\'e}briquement clos. On peut
enlever de $\sX$ les composantes connexes de $\sX_k$ qui ne
rencontrent pas $g(S)$, ferm{\'e}es dans $\sX$, et donc supposer
$\sX_k$ connexe. On peut ensuite remplacer $\sX$ par sa
composante connexe contenant $g(S)$ et donc supposer $\sX$
connexe. Il suit alors du
lemme~\ref{lemm:connexite_et_fibre_generique} que $\sX_K$ est
connexe, donc $C=|\sX|$ est ouvert, ce qui termine la preuve.
\end{adoc}

\subsection{Composantes irr{\'e}ductibles le long d'une
section~: preuve}
\label{ss_section:ci_preuve}

Dans cette section, nous d\'emontrons la
proposition~\ref{comp_irred_d_une_section}.

\begin{adoc}{Remarques pr{\'e}liminaires} \label{adoc:remarques}
Il est naturel d'essayer de prouver la
proposition~\ref{comp_irred_d_une_section} en deux {\'e}tapes,
d'abord dans le cas o{\`u} $\sX$ est un sch{\'e}ma, puis en
ramenant le cas g{\'e}n{\'e}ral {\`a} ce cas en utilisant un atlas
et diverses techniques semblables {\`a} celles
de~\ref{preuve_de_Brochard}. Ces m{\'e}thodes font jouer un r{\^o}le
important {\`a} des parties $X^0(e)$, r{\'e}unions de composantes
connexes de $X_s$ qui rencontrent $e(V)$, o{\`u} $e:V\to \sX$ est
un $S$-morphisme depuis un $S$-sch{\'e}ma $V$ distinct de $S$,
typiquement un morphisme lisse surjectif issu d'un atlas de
$\sX$. On est amen{\'e} {\`a} consid{\'e}rer $X^1(e)$, r{\'e}union
des composantes irr{\'e}ductibles ouvertes de $X_s$ qui
rencontrent $e(V)$, pour des $e:V\to \sX$ {\`a} valeurs dans
le lieu unicomposante. On rencontre alors les probl{\`e}mes
suivants~:

\begin{trivlist}
\itemn{i} On ne peut pas toujours relever des points du lieu
unicomposante en des points du lieu unicomposante, m{\^e}me
par des morphismes finis {\'e}tales. Un exemple est donn{\'e} par
le rev{\^e}tement double non trivial de la cubique nodale
(voir figure~\ref{fig1}).

\begin{figure}[ht]
\centerline{\begin{tikzpicture}
\draw[thick] (0,2) .. controls (0,4) and (2,4) .. (4.5,0.5);
\draw[thick] (4.5,3.5) .. controls (2.5,0) and (0.3,0) .. (0.3,2);
\draw[draw=white,double=black,ultra thick]
(0.3,2) .. controls (0.3,4) and (2.5,3) .. (4,0.5);
\draw[draw=white,double=black,ultra thick]
(4,3.5) .. controls (2,0) and (0,0) .. (0,2);
\filldraw [black] (2.92,1.92) circle (2pt);
\filldraw [black] (3.38,1.9) circle (2pt);
\draw [->,thick] (6,2) -- (8,2);
\draw[thick] (13,3.5) .. controls (11,0) and (9,0) .. (9,2);
\draw[thick] (9,2) .. controls (9,4) and (11,4) .. (13,0.5);
\filldraw [black] (11.98,1.98) circle (2pt);
\end{tikzpicture}}
\caption{{\it Rev{\^e}tement double de la cubique nodale}}
\label{fig1}
\end{figure}
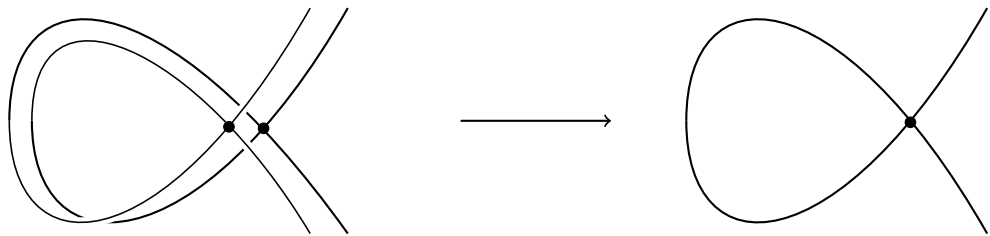

\itemn{ii} La formation de $X^1(e)$ ne commute pas au changement
de base. Ceci est li{\'e} au fait que la formation des
composantes irr{\'e}ductibles ouvertes sans point rationnel
ne commute pas au changement de base. Ce ph{\'e}nom{\`e}ne se
produit m{\^e}me si $V\to S$ est fini {\'e}tale~: prendre
pour $X$ la conique sur $S=\Spec(\bQ)$ d'{\'e}quation
$x^2-2y^2=0$, et $e:\Spec(\bQ(\sqrt{2})\to X$ donn{\'e}
par le point $(x,y)=(\sqrt{2},1)$.
\end{trivlist}
{\`A} cause de ces difficult{\'e}s, nous ne pr{\'e}sentons pas
de preuve de la proposition~\ref{comp_irred_d_une_section}
analogue {\`a}~\ref{preuve_de_Brochard}.
\end{adoc}

\begin{adoc}{Preuve de~\ref{comp_irred_d_une_section}}
Nous d{\'e}montrons maintenant~\ref{comp_irred_d_une_section}
en suivant une strat{\'e}gie classique, comme
dans~\ref{preuve_classique}. Le fait que la formation du lieu
unicomposante relatif commute au changement de base est une
cons{\'e}quence directe de sa d{\'e}finition.
Pour v{\'e}rifier la propri{\'e}t{\'e} analogue pour la composante
irr{\'e}ductible ouverte le long d'une section, nous utiliserons le
lemme suivant, qui est une variante de~\cite{EGA}~IV.4.5.13.

\begin{lemm} \label{lemm:geom_integre}
Soient $k$ un corps et $f:\sY\to\sX$ un morphisme de
$k$-champs alg{\'e}briques. On suppose que $\sX$ est
irr{\'e}ductible et ne poss{\`e}de qu'un nombre fini de
composantes irr{\'e}ductibles g{\'e}om{\'e}tri\-ques, que $\sY$ est
g{\'e}om{\'e}triquement irr{\'e}ductible et que la pr{\'e}image par
$f$ du lieu unicomposante de $\sX/k$ est non vide.
Alors, $\sX$ est g{\'e}om{\'e}triquement irr{\'e}ductible.
\end{lemm}

\begin{demo}
Soit $k'$ une cl{\^o}ture alg{\'e}brique de $k$~; on note avec
un <<~${}'$~>> toutes les donn{\'e}es obtenues par changement de
base de $k$ {\`a} $k'$, dont le morphisme $f':\sY'\to\sX'$.
Soit $p:\sX'\to\sX$ la projection, qui est
ouverte et ferm{\'e}e. Obervons que le lieu unicomposante $\sU$
de $\sX/k$ est un sous-champ ouvert, puisque c'est l'image
par~$p$ du lieu unicomposante de $\sX'/k'$ qui est ouvert.
La d{\'e}monstration ne fait intervenir que les espaces topologiques
sous-jacents {\`a} $\sX,\sY,\sU$ et {\`a} leurs homologues sur
$k'$~; on les note $X,Y,U$, etc. Si $X'$ n'est pas irr{\'e}ductible,
on peut partager ses composantes irr{\'e}ductibles en deux paquets
disjoints et former les ferm{\'e}s $F',G'$ r{\'e}unions des
composantes irr{\'e}ductibles de chacun des deux paquets, et les
ouverts disjoints $A'=X'\setminus G'$, $B'=X'\setminus F'$. On
note $F,G,A,B$ les images dans~$X$. Puisque~$p$ est ouverte et
ferm{\'e}e, on voit que $A$ et $B$ sont des ouverts denses de $X$
et que $F$ et $G$ sont des ferm{\'e}s d'int{\'e}rieur non vide, donc
$F=G=X$. De plus, comme $p(X'\setminus U')\subset X\setminus U$
on trouve $U\subset A$ et $U\subset B$. Puisque par hypoth{\`e}se
$f^{-1}(U)$ est non vide, il d{\'e}coule de ce qui pr{\'e}c{\`e}de
que les ouverts disjoints $(f')^{-1}(A')$ et $(f')^{-1}(B')$
sont non vides. Ceci n'est pas possible, car $Y'$ est suppos{\'e}
irr{\'e}ductible.
\end{demo}

Soient $S'\to S$ un morphisme, $\sX'=\sX\times_S S'$,
$g':S'\to\sX'$ la section d{\'e}duite de la section $g$,
et~$I'$ la r{\'e}union des $|\sI_{s'}(g')|$. V{\'e}rifions que
la partie $I'$ est la pr{\'e}image de $I$ par l'application
continue $|\sX'|\to |\sX|$. On le v{\'e}rifie fibre {\`a} fibre.
Si $s'\in S'$ a pour image $s\in S$, la composante irr{\'e}ductible
ouverte $\sI_s$ est g{\'e}om{\'e}triquement irr{\'e}ductible
d'apr{\`e}s le lemme~\ref{lemm:geom_connexe} appliqu{\'e} {\`a} la
section $g$. Utilisant la caract{\'e}risation des composantes
irr{\'e}ductibles ouvertes comme ouverts irr{\'e}ductibles maximaux
et le fait que le morphisme $\sX'_{s'}\to \sX_s$ est ouvert, on
voit que $\sI'_{s'}=\sI_s\otimes_{k(s)}k(s')$. L'assertion
{\`a} d{\'e}montrer en d{\'e}coule.

Il nous reste {\`a} d{\'e}montrer que $U$ et $I$ sont ouverts,
ce que nous ferons en m{\^e}me temps. Comme la formation de $U$ et
$I$ commute au changement de base, on peut travailler
{\'e}tale-localement sur $S$. En particulier, on peut supposer
que $S=\Spec(A)$ est affine. Comme $\sX\to S$ est de pr{\'e}sentation
finie, on peut ensuite supposer que $A$ est noeth{\'e}rien.

Montrons que $U$ et $I$ sont localement constructibles.
On se ram{\`e}ne imm{\'e}diatement au cas o{\`u}~$S$ est int{\`e}gre de
point g{\'e}n{\'e}rique $\eta$. Soient $\sZ_{0,\eta},\dots,\sZ_{d,\eta}$
les composantes irr{\'e}ductibles de $\sX_\eta$, en supposant que
$g(\eta)\in \sZ_{0,\eta}$ dans le cas ii). Quitte {\`a} faire une
extension {\'e}tale de $S$, on peut supposer que $\sZ_{i,\eta}$ est
g{\'e}om{\'e}triquement int{\`e}gre pour tout $i$. Soit $\sZ_i$ un
sous-champ ferm{\'e} de $\sX$ dont la fibre au point $\eta$ est
$\sZ_{i,\eta}$. D'apr{\`e}s
\ref{prop:EGA_95}~(i), (ii), (iv) et \ref{theo:EGA_977}~(iv), quitte
{\`a} remplacer $S$ par un voisinage de $\eta$, on peut supposer que
les $\sZ_i$ recouvrent $\sX$ et que pour tout $s\in S$, leurs fibres
sont des ferm{\'e}s g{\'e}om{\'e}triquement int{\`e}gres tels que
$\sZ_{i,s}\not\subset \sZ_{j,s}$ pour $i\ne j$. Par ailleurs dans le cas
ii) on a $g^{-1}(\sZ_0)=S$ car ce ferm{\'e} contient le point
g{\'e}n{\'e}rique de $S$ qui est irr{\'e}ductible, et $g(S)$ est inclus
dans le compl{\'e}mentaire dans $\sX$ de $\cup_{i\ne 0}\sZ_i$. Il est
alors clair que $U$ est {\'e}gal au
compl{\'e}mentaire de $\cup_{i,j}|\sZ_i\cap \sZ_j|$ et que $I$ est
{\'e}gal au compl{\'e}mentaire de $\cup_{i\ne 0}|\sZ_i|$, qui sont des
parties ouvertes, donc localement constructibles.

On conclut enfin que $U$ et $I$ sont ouverts. D'apr{\`e}s ce qui
pr{\'e}c{\`e}de, on peut supposer que~$S$ est le spectre d'un anneau
de valuation discr{\`e}te $(R,K,k,\pi)$ que l'on peut supposer
complet {\`a} corps r{\'e}siduel alg{\'e}briquement clos.
On peut aussi enlever de $\sX$ ses composantes irr{\'e}ductibles
incluses dans $\sX_k$, puisqu'elles sont ferm{\'e}es et ne
rencontrent pas $\sI_k$, et donc supposer $\sX$ plat sur $S$.
Il suffit de montrer que $U$ et $I$ sont ouverts au voisinage
d'une quelconque composante irr{\'e}ductible $\sY$ de $\sX_k$, que
l'on choisit comme {\'e}tant la composante contenant
$g(\Spec(k))$ dans le cas ii). On peut donc enlever de $\sX$
toutes les composantes irr{\'e}ductibles de $\sX_k$ distinctes de
$\sY$ et supposer $\sX_k$ g{\'e}om{\'e}triquement int{\`e}gre.
On peut ensuite enlever les composantes irr{\'e}ductibles de
$\sX$ incluses dans $\sX_K$ et donc supposer $\sX$ pur sur $S$
(pour des d{\'e}tails sur la puret{\'e}, voir
l'annexe~\ref{annexe:purete}, \cite{RG} ou \cite{Ro2}). Il suit
alors du th.~\ref{theo:CC_CI_pour_les_champs}~(iii) que $\sX_K$
est int{\`e}gre, de sorte que $U=I=|\sX|$ et notre assertion
est prouv{\'e}e.
\end{adoc}

\subsection{Foncteurs des composantes ouvertes}
\label{ss_section:representabilite_de_pi_0}

Nous allons d{\'e}montrer que pour les champs plats, de pr{\'e}sentation
finie et {\`a} fibres g{\'e}om{\'e}triquement r{\'e}duites, les foncteurs
$\pi_0(\sX/S)$ et $\Irr(\sX/S)$ sont repr{\'e}sentables par des
espaces alg{\'e}briques {\'e}tales et quasi-compacts. Nous pr{\'e}cisons
tout d'abord un petit point de th{\'e}orie
(\ref{adoc:quotient_par_un_groupoide}) avant d'{\'e}noncer
le th{\'e}or{\`e}me (\ref{theo:pi_0_representable}).

\begin{adoc}{Quotient d'un champ alg{\'e}brique par un groupo{\"\i}de}
\label{adoc:quotient_par_un_groupoide}
Nous aurons besoin ci-dessous de consid{\'e}rer le quotient d'un
champ alg{\'e}brique par une relation d'{\'e}quivalence {\'e}tale.
Cette op{\'e}ration ne pr{\'e}sente pas v{\'e}ritablement de
difficult{\'e}, et les sp{\'e}cialistes de l'homotopie savent
m{\^e}me faire beaucoup mieux (voir notamment \cite{TV},
paragraphe~1.3.4). Malheureusement, leurs travaux restent
difficiles \`a lire pour de nombreux g\'eom\`etres, y compris par
le r\'edacteur de ces lignes. Faute de r{\'e}f{\'e}rence dans
laquelle ces r{\'e}sultats sont {\'e}nonc{\'e}s dans un cadre
plus simple, adapt{\'e} {\`a} nos besoins, nous esquisserons ces
constructions.

Soit~$S$ un sch{\'e}ma. Nous appellerons {\em groupo{\"\i}de plat} (en champs
alg{\'e}briques) et nous noterons $\sX_1\rightrightarrows \sX_0$ la
donn{\'e}e de deux $S$-champs alg{\'e}briques $\sX_0,\sX_1$ et de~:
\begin{trivlist}
\itemm{1} deux morphismes fid{\`e}lement plats et de pr{\'e}sentation
finie $s,t:\sX_1\to\sX_0$,
\itemm{2} un morphisme de composition $c:\sX_1\times_{t,\sX_0,s}\sX_1\to\sX_1$,
\itemm{3} un morphisme d'identit{\'e} $e:\sX_0\to\sX_1$,
\itemm{4} un morphisme d'inversion $i:\sX_1\to\sX_1$,
\end{trivlist}
ainsi que d'un certain nombre de $2$-isomorphismes de
compatibilit{\'e} exprimant la $2$-commutativit{\'e} des diagrammes
des axiomes bien connus de groupo{\"\i}de (voir par exemple
\cite{LMB}, 2.4.3). Par souci de l{\'e}g{\`e}ret{\'e}, nous ne pr{\'e}cisons
pas ces $2$-isomorphismes, mais ils peuvent facilement {\^e}tre
{\'e}crits en s'inspirant par exemple de \cite{Ro1}, section~1.
Comme dans la definition~1.3 de~\cite{Ro1}, on a une notion de
groupo{\"\i}de {\em strict} correspondant au cas o{\`u} tous ces
$2$-isomorphismes sont l'identit{\'e}. \'Etant donn\'e un champ
$\sY$, un {\em morphisme $\sX_1$-invariant de $\sX_0$ vers $\sY$}
consiste en la donn\'ee d'un morphisme $f:\sX_0\to\sY$ et d'un
$2$-isomorphisme $\sigma:f\circ s\Rightarrow f\circ t$ compatible
avec les $2$-isomorphismes de structure du groupo\"{i}de.
On peut former le {\em stabilisateur} $\sS$ d{\'e}fini par le
$2$-produit fibr{\'e}
$$
\xymatrix{
\sS \ar[r] \ar[d] & \sX_0 \ar[d]^\Delta \\
\sX_1 \ar[r]^>>>>>{s,t} & \sX_0\times_S\sX_0 \\}
$$
o{\`u} $\Delta$ est la diagonale. L'{\'e}nonc{\'e} de repr{\'e}sentabilit{\'e}
qui nous int{\'e}resse est que si $\sS\to\sX_0$ est un morphisme
repr{\'e}sentable, il existe un {\em champ alg{\'e}brique quotient}
$\pi:\sX_0\to\sX:=\sX_0/\sX_1$ qui v{\'e}rifie la $2$-propri{\'e}t{\'e}
universelle suivante~: pour tout morphisme de champs alg{\'e}briques
$f:\sX_0\to\sY$ qui est $\sX_1$-invariant, il existe un morphisme
$f':\sX\to\sY$ tel que $f=f'\circ\pi$, unique {\`a} un unique
$2$-isomorphisme pr{\`e}s.

On peut d{\'e}crire $\sX_0/\sX_1$ de la mani{\`e}re suivante~:
c'est le champ associ{\'e} {\`a} la cat{\'e}gorie fibr{\'e}e en
groupo{\"\i}des $\sP$ dont les fibres sont les cat{\'e}gories
$\sP(T)$, pour des $S$-sch{\'e}ma~$T$ variables, telles
que $\Obj(\sP(T))=\sX_0(T)$ et $\Hom_{\sP(T)}(x,y)$ est
l'ensemble des paires $(x_1,\varphi)$ o{\`u} $x_1\in\sX_1(T)$
v{\'e}rifie $s(x_1)=x$ et $\varphi:t(x_1)\to y$ est un morphisme
dans $\sX_0(T)$.
La composition des morphismes dans $\sP(T)$ se fait ainsi~:
{\'e}tant donn{\'e}s un morphisme $(x_1,\varphi)$ entre $x$ et $y$
et un morphisme $(y_1,\psi)$ entre $y$ et $z$, le triplet
$(x_1,y_1,\varphi)$ est un objet de $\sX_1\times_{t,\sX_0,s}\sX_1$
et on peut poser $x'_1:=c(x_1,y_1,\varphi)$. Le compos{\'e} de
$(x_1,\varphi)$ et $(y_1,\psi)$ est alors $(x'_1,\psi)$.

Le r\'esultat de quotient que nous venons d'esquisser recouvre
le quotient d'un champ alg\'ebrique par l'action d'un sch\'ema
en groupes comme dans \cite{Ro1}, par l'action libre d'un champ
alg\'ebrique en groupes comme dans \cite{La}, ou par une relation
d'\'equivalence plate, comme ci-dessous.
\end{adoc}

\begin{theo} \label{theo:pi_0_representable}
Soit $\sX$ un $S$-champ alg{\'e}brique plat, de pr{\'e}sentation
finie, {\`a} fibres g{\'e}om{\'e}tri\-quement r{\'e}duites.
\begin{trivlist}
\itemn{i} Les foncteurs $\pi_0(\sX/S)$ et $\Irr(\sX/S)$
sont repr{\'e}sentables par des espaces alg{\'e}briques {\'e}tales
et quasi-compacts sur~$S$.
\itemn{ii} Soit la relation d'{\'e}quivalence d{\'e}finie comme
sous-cat{\'e}gorie pleine $\sR\subset\sX\times\sX$ telle que
deux points $u,v:T\to\sX$ sont {\'e}quivalents ssi
pour tout point g{\'e}om{\'e}trique $t:\Spec(\Omega)\to T$,
les points $u(t)$ et $v(t)$ sont dans la m{\^e}me composante connexe
de $\sX_\Omega$. Cette relation est repr{\'e}sentable par la
\cco\ de $\sX\times\sX$ le long de la section diagonale.
De plus, il existe un morphisme $\sX\to \pi_0(\sX/S)$ gr{\^a}ce
auquel $\sX$ s'identifie {\`a} la \cco\ universelle et
$\pi_0(\sX/S)$ au quotient $\sX/\sR$.
\itemn{iii} Soit $\sU\subset\sX$ le lieu unicomposante de
$\sX$ sur $S$, et soit la relation d'{\'e}quivalence
$\sS\subset\sU\times\sU$, sous-cat{\'e}gorie pleine telle que
deux points $u,v:T\to\sU$ sont {\'e}quivalents ssi
pour tout point g{\'e}om{\'e}trique $t:\Spec(\Omega)\to T$,
les points $u(t)$ et $v(t)$ sont dans la m{\^e}me composante
irr{\'e}ductible ouverte de $\sU_\Omega$. Cette relation est
repr{\'e}sentable par la \cio\ de $\sU\times\sU$ le long de la
section diagonale.
De plus, il existe un morphisme $\sU\to \Irr(\sX/S)$ gr{\^a}ce
auquel $\sU$ s'identifie {\`a} la \cio\ universelle et
$\Irr(\sX/S)$ au quotient $\sU/\sS$.
\end{trivlist}
\end{theo}

Rappelons que le lieu unicomposante $\sU$ est un sous-champ
ouvert de $\sX$, voir la proposition~\ref{comp_irred_d_une_section}.

\begin{demo}
(i) On utilise les crit{\`e}res de repr{\'e}sentabilit{\'e}
d'Artin (corollaire~5.2 de \cite{Ar1}) pour les foncteurs
$F=\pi_0(\sX/S)$ et
$F=\Irr(\sX/S)$. Le raisonnement {\'e}tant le m{\^e}me pour les
deux foncteurs consid{\'e}r{\'e}s, disons que $F=\pi_0(\sX/S)$
pour fixer les id{\'e}es. On se ram{\`e}ne par les arguments
habituels
au cas o{\`u}~$S$ est le spectre d'une
$\bZ$-alg{\`e}bre de type fini. D'apr{\`e}s le
lemme~\ref{lemme:faisceau_etale}, le foncteur $F$ est un
faisceau localement de pr{\'e}sentation finie {\`a} diagonale
ouverte (condition $(1)$ d'Artin) et de plus formellement
{\'e}tale sur~$S$, donc les th{\'e}ories de d{\'e}formations
et d'obstructions sont nulles (conditions $(2)$ et $(4)$ d'Artin).

Il ne reste que la condition~$(3)$ d'effectivit{\'e} des composantes
formelles {\`a} v{\'e}rifier~: on doit montrer que pour tout anneau
local complet $(R,m)$ et tout morphisme $\Spec(R)\to S$ qui induit
une extension de corps r{\'e}siduels de type fini, l'application
$F(R)\to\varprojlim F(R/m^n)$ est injective et d'image dense.
Puisque les applications $F(R/m^{n+1})\to F(R/m^n)$ sont des
bijections, il s'agit juste de montrer que l'application
$\varphi:F(R)\to F(k)$ est bijective.

Soient $\sZ,\sZ'$ deux \cco\ de $\sX_R$ sur $R$ telles
que $\sZ_k=\sZ'_k$. Il suffit de montrer que $\sZ=\sZ'$ apr{\`e}s
une extension {\'e}tale surjective de $R$. Comme $\sX\to S$
est {\`a} fibres g{\'e}om{\'e}triques r{\'e}duites, son lieu de
lissit{\'e} est un ouvert $S$-sch{\'e}matiquement dense.
Par~\cite{EGA}~IV.17.16.3, il existe une extension finie
{\'e}tale locale $(R',m',k')$ de $R$ et une section
$g:\Spec(R')\to \sX_{R'}$ qui envoie le point ferm{\'e}
dans $\sZ_k=\sZ'_k$. Il est alors clair que $\sZ$ et $\sZ'$
sont \'egales \`a $\sC(g)$, la \cco\ de $\sX_{R'}$ le long
de $g$.

Enfin, soit $\sZ$ une \cco\ de $\sX_k$ sur $k$. On consid{\`e}re
de nouveau une extension finie {\'e}tale locale $(R',m',k')$
de $R$ et une section $g:\Spec(R')\to \sX_{R'}$ qui envoie le
point ferm{\'e} dans $\sZ$. Il suit alors du
lemme~\ref{comp_connexe_d_une_section} que la composante
connexe $\sC(g)$ le long de $g$ est une \cco\ de
$\sX_{R'}$ qui rel{\`e}ve $\sZ_{k'}$.
D'apr{\`e}s le lemme~\ref{lemm:geom_connexe}, son image
dans $\sX_R$ est une \cco\ de $\sX_R$ qui rel{\`e}ve $\sZ$, donc
$\varphi$ est surjective.

\medskip

\noindent (ii)
{\'E}tant donn{\'e} que pour tout point g{\'e}om{\'e}trique $\Spec(\Omega)\to S$
les composantes connexes de $\sX_\Omega\times\sX_\Omega$ sont les produits
$\sC_i\times\sC_j$ des composantes connexes de $\sX_\Omega$, il est
imm{\'e}diat que la relation $\sR$ est repr{\'e}sentable par la \cco\ de
$\sX\times\sX$ le long de la section diagonale. Par ailleurs, en
associant {\`a} une section de $\sX$ au-dessus d'un sch{\'e}ma $T$
la \cco\ de $\sX\times_S T$ le long de cette section, on d{\'e}finit
un morphisme surjectif $a:\sX\to \pi_0(\sX/S)$. Comme $\sX\to S$ est
plat et $\pi_0(\sX/S)\to S$ {\'e}tale, le morphisme $a$ est plat
et donc un {\'e}pimorphisme fppf. Il est clair que $a$ passe au quotient
en un monomorphisme $\sX/\sR\to \pi_0(\sX/S)$ qui est donc
une immersion ouverte puis un isomorphisme.

\medskip

\noindent (iii) Une fois qu'on a remarqu{\'e} que le lieu
unicomposante de $\sX\times\sX$ est $\sU\times\sU$, la preuve
est formellement la m{\^e}me que celle du point (ii) en
utilisant la proposition~\ref{comp_irred_d_une_section}
au lieu de la proposition~\ref{comp_connexe_d_une_section}.
\end{demo}

\begin{remas}
(1) Il est facile de voir que $\sX$ est un $S$-espace alg\'ebrique
\'etale et quasi-compact si et seulement si $\sX\to\pi_0(\sX/S)$
est un isomorphisme. En particulier, on peut obtenir ainsi un
exemple dans lequel $\pi_0(\sX/S)$ n'est pas un sch\'ema.

\no (2) Si $\pi_0(\sX/S)$ est s{\'e}par{\'e}, alors c'est un
sch{\'e}ma par~\cite{Kn}, chap.~II, cor.~6.17. Mais ce foncteur
est rarement s{\'e}par{\'e}. Par exemple, dans les conditions
du th\'eor\`eme, si $S$ est un trait strictement hens\'elien et
$\sX$ est connexe, il est facile de voir que $\pi_0(\sX/S)$
est le sch\'ema obtenu en recollant~$n$ copies de $S$ le long
de leur point ouvert commun, o\`u $n$ est le nombre de composantes
connexes de la fibre sp\'eciale de $\sX$. Ainsi $\pi_0(\sX/S)$
n'est pas s\'epar\'e si $n\ge 2$.

\no (3) Il y a une description birationnelle de $\Irr(\sX/S)$
qui m{\'e}rite d'{\^e}tre mentionn{\'e}e. Soit le foncteur $F$ dont
les valeurs, pour chaque $S$-sch{\'e}ma $T$, sont les classes
d'{\'e}quivalence de sous-champs ouverts $\sU\subset\sX_T$,
fid{\`e}lement plats et de pr{\'e}sentation finie sur $T$, {\`a}
fibres g{\'e}om{\'e}triquement irr{\'e}ductibles. Deux tels ouverts
$\sU,\sV$ sont {\'e}quivalents si et seulement si $\sU\cap\sV$
est surjectif sur $T$. Lorsque $\sX$ est plat, de pr{\'e}sentation
finie, {\`a} fibres g{\'e}om{\'e}triquement r{\'e}duites,
il est clair que $F$ est isomorphe {\`a} $\Irr(\sX/S)$~; {\`a}
l'aide de la proposition~\ref{comp_irred_d_une_section},
on d{\'e}finit deux morphismes en sens oppos{\'e}s, inverses l'un
de l'autre.
\end{remas}

\subsection{Fonctorialit{\'e}} \label{subsection:fonctorialite}

\begin{coro}
Soit $\sX$ un $S$-champ alg{\'e}brique plat, de pr{\'e}sentation
finie, {\`a} fibres g{\'e}om{\'e}tri\-quement r{\'e}duites. Alors
il existe un morphisme surjectif $\Irr(\sX/S)\to\pi_0(\sX/S)$.
\end{coro}

\begin{demo}
Il est clair que le morphisme compos{\'e} $\sU\into\sX\to\pi_0(\sX/S)$
passe au quotient par la relation d'{\'e}quivalence $\sS$.
\end{demo}

Nous appellerons {\em application $S$-rationnelle} de $\sX$ vers $\sY$
une classe d'{\'e}quivalence de morphismes $\sU\to\sY$ d{\'e}finis sur des
ouverts $S$-denses de $\sX$, o{\`u} l'on convient que deux morphismes
$\sU\to\sY$ et $\sV\to\sY$ sont {\'e}quivalents s'ils co{\"\i}ncident
sur $\sU\cap \sV$. On trouve dans~\cite{EGA}~IV.20 la terminologie
de {\em pseudo-morphisme de $\sX$ dans $\sY$ relativement {\`a} $S$}.

\begin{coro} \label{coro:covariance_en_X}
Soit $f:\sX\dashrightarrow\sY$ une application $S$-rationnelle
entre $S$-champs alg{\'e}briques plats, de pr{\'e}sentation finie,
{\`a} fibres g{\'e}om{\'e}triquement r{\'e}duites.
\begin{trivlist}
\itemn{i}
$f$ induit un morphisme $\pi_0(f):\pi_0(\sX/S)\to\pi_0(\sY/S)$.
\itemn{ii}
Si $f^{-1}(\sU_\sY)\subset\sU_\sX$, $f$ induit un
morphisme $\Irr(f):\Irr(\sX/S)\to\Irr(\sY/S)$.
\end{trivlist}
\end{coro}

On notera que, compte tenu du caract{\`e}re birationnel de
$\Irr(\sX/S)$, l'assertion de (ii) est encore valable sous
l'hypoth{\`e}se plus faible que $f^{-1}(\sU_\sY)\cap\sU_\sX$ est
$S$-dense dans $\sU_\sX$, par restriction {\`a} cet ouvert.

\begin{demo}
Pour un champ $\sX$, notons $\sU_\sX$, $\sR_\sX$, $\sS_\sX$
les objets d{\'e}finis dans le
th{\'e}or{\`e}me~\ref{theo:pi_0_representable}.
Supposons d'abord que $f$ est partout d{\'e}finie, i.e. un
morphisme $f:\sX\to\sY$. Dans ce cas, il est clair que le
morphisme compos{\'e} $\sX\to\sY\to\pi_0(\sY/S)$ passe au quotient
par la relation d'{\'e}quivalence $\sR_\sX$ en un morphisme
$\pi_0(\sX/S)\to\pi_0(\sY/S)$. De m{\^e}me, si
$f^{-1}(\sU_\sY)\subset\sU_\sX$ alors le morphisme compos{\'e}
$\sU_\sX\to\sU_\sY\to\Irr(\sY/S)$ passe au quotient par la
relation d'{\'e}quivalence $\sS_\sX$ en un morphisme
$\Irr(\sX/S)\to\Irr(\sY/S)$.

Si $f$ est une immersion ouverte $S$-dense, les morphismes
ci-dessus sont {\'e}tales et bijectifs sur les fibres
g{\'e}om{\'e}triques, donc ce sont des isomorphismes. On en
d{\'e}duit que tout ce qui vient d'{\^e}tre dit est encore
valable pour une application $S$-rationnelle
$f:\sX\dashrightarrow\sY$ quelconque, ce qui termine de
d{\'e}montrer (i) et (ii).
\end{demo}

\begin{coro}
Soit $\sX$ un $S$-champ alg{\'e}brique plat, de pr{\'e}sentation
finie, {\`a} fibres g{\'e}om{\'e}tri\-quement r{\'e}duites. On
suppose que $\sX$ est un champ alg{\'e}brique mod{\'e}r{\'e} au sens
de {\rm \cite{AOV}} et on note $p:\sX\to X$ son espace grossier de
modules. Alors $X$ est un $S$-espace alg{\'e}brique plat, de
pr{\'e}sentation finie, {\`a} fibres g{\'e}om{\'e}triquement
r{\'e}duites et l'application induite $\pi_0(\sX/S)\to\pi_0(X/S)$
est un isomorphisme.
\end{coro}

\begin{demo}
D'apr{\`e}s \cite{AOV}, corollary 3.3, l'espace $X$ est plat sur $S$
et ses fibres g{\'e}om{\'e}triques sont les espaces de modules des
fibres g{\'e}om{\'e}triques de $\sX$. Ceci montre que $X$ est {\`a}
fibres g{\'e}om{\'e}triquement r{\'e}duites. L'application induite
par fonctorialit{\'e} de $p:\sX\to X$ est clairement bijective sur
les fibres g{\'e}om{\'e}triques, donc c'est un isomorphisme.
\end{demo}

\section{Composantes ferm{\'e}es}

\subsection{D{\'e}finitions}

Pour un champ $\sX\to S$ {\`a} fibres non (g{\'e}om{\'e}triquement)
r{\'e}duites, le foncteur des \cco\ n'est pas n{\'e}cessairement
repr{\'e}sentable. Par exemple, consid{\'e}rons un anneau de valuation
discr{\`e}te complet {\`a} corps r{\'e}siduel alg{\'e}briquement clos
$(R,K,k,\pi)$ et prenons pour $\sX$ la r{\'e}union de deux sections de
la droite affine sur $R$ qui se rencontrent dans la fibre sp{\'e}ciale,
par exemple le spectre de $R[x]/(x^2-\pi^n x)$. Il est clair que
$\pi_0(\sX/S)(R)=\emptyset$ de sorte que l'approximation des
composantes connexes formelles est prise en d{\'e}faut. Il en
d{\'e}coule que $\pi_0(\sX/S)$ n'est pas repr{\'e}sentable par
un espace alg{\'e}brique. N{\'e}anmoins, dans cet exemple le champ
$\sX$ est propre. Gr{\^a}ce {\`a} l'existence du sch{\'e}ma de
Hilbert, pour ces champs, on a des notions utiles de composantes
connexes et irr{\'e}ductibles {\em ferm{\'e}es}, ainsi que nous
le voyons maintenant.

\begin{defis} \label{defi:ccrf_cirf}
Soit $\sX$ un champ alg{\'e}brique de pr{\'e}sentation finie sur un
sch{\'e}ma $S$.
\begin{trivlist}
\itemn{1} Une {\em composante connexe (relative) ferm{\'e}e
(en abr{\'e}g{\'e} \ccf) de $\sX$ sur $S$} est un sous-champ
$\sC\subset\sX$ ferm{\'e}, plat et de pr{\'e}sentation finie sur~$S$,
tel que le support de $\sC_s$ est une composante connexe de $\sX_s$,
pour tout point g{\'e}om{\'e}trique $s$ de $S$. On note
$\pi_0(\sX/S)^\ff$ le foncteur des \ccf\ de $\sX$ sur $S$.
\itemn{2} On dit qu'une \ccf\ $\sC$ est {\em r{\'e}duite}
si ses fibres g{\'e}om{\'e}triques sont r{\'e}duites.
On note $\pi_0(\sX/S)^\fr$
le sous-foncteur correspondant de $\pi_0(\sX/S)^\ff$.
\itemn{3} Une {\em composante irr{\'e}ductible (relative) ferm{\'e}e
(en abr{\'e}g{\'e} \cif) de $\sX$ sur $S$} est un sous-champ
$\sI\subset\sX$ ferm{\'e}, plat et de pr{\'e}sentation finie sur~$S$,
tel que le support de $\sI_s$ est une composante
irr{\'e}ductible de $\sX_s$, pour tout point g{\'e}om{\'e}trique
$s$ de $S$. On note $\Irr(\sX/S)^\ff$ le foncteur des \cif\
de $\sX$ sur $S$.
\end{trivlist}
\end{defis}

Lorsque $\sX/S$ est {\`a} fibres g{\'e}om{\'e}triquement r{\'e}duites,
ces foncteurs se comparent aux foncteurs de composantes ouvertes
(voir proposition~\ref{prop:comparaison_ouverte_fermee}).
Dans le cas contraire, on s'attend bien s{\^u}r
{\`a} ce que les foncteurs de composantes ferm{\'e}es soient
ramifi{\'e}s. L'exemple~\ref{exem:non_representabilite}(1) ci-dessous
montre que $\pi_0(\sX/S)^\ff$ peut {\^e}tre d{\'e}mesur{\'e}ment
gros, de dimension positive et non quasi-compact, m{\^e}me
si $\sX/S$ est propre. Ceci est d{\^u} {\`a} la pr{\'e}sence de
composantes immerg{\'e}es dans les \ccf\
C'est ce qui justifie l'introduction de $\pi_0(\sX/S)^\fr$
qui sera, lui, quasi-fini. Pour les composantes irr{\'e}ductibles,
la situation est diff{\'e}rente car,
comme on va le voir, en pr{\'e}sence de fibres non r{\'e}duites,
on ne sait pas en fabriquer de bons espaces de modules.

\begin{exems} \label{exem:non_representabilite}
(1) Soient $k$ un corps et $X_0$ un $k$-sch{\'e}ma
g{\'e}om{\'e}triquement connexe. Soit $X=X_0[\epsilon]$
obtenu {\`a} partir de $X$ par le changement de base de $k$
{\`a} l'anneau des nombres duaux $k[\epsilon]/(\epsilon^2)$.
On a alors un isomorphisme de foncteurs
$\Hilb(X_0/k)\simeq\pi_0(X/k)^\ff$
qui envoie le sous-sch{\'e}ma ferm{\'e} d'id{\'e}al $\cI$
sur la \ccf\ d'id{\'e}al $\epsilon\cI$.

\no (2) Si $\sX/S$ n'est pas s{\'e}par{\'e}, les {\'e}l{\'e}ments
formels des foncteurs $\pi_0(\sX/S)^\ff$ et $\Irr(\sX/S)^\ff$
ne peuvent pas {\^e}tre approxim{\'e}s en g{\'e}n{\'e}ral, m{\^e}me avec $\sX/S$ plat et
pur. Par exemple, prenons pour base un anneau de valuation discr{\`e}te
$(R,K,k,\pi)$ complet {\`a} corps r{\'e}siduel alg{\'e}briquement clos et prenons
pour $X$ le sch{\'e}ma obtenu en recollant deux copies $U,V$ de la droite
affine sur~$R$ le long de leurs fibres g{\'e}n{\'e}riques.
Pour tout $n\ge 1$, posons $R_n=R/\pi^n$, $X_n=X\otimes R_n$,
$U_n=U\otimes R_n$, $V_n=V\otimes R_n$. Alors les ouverts
$U_n$ et $V_n$ sont disjoints dans $X_n$, et la collection $(U_n)$
d{\'e}finit une \ccf\ (ou une \cif) {\'e}l{\'e}ment de
$\varprojlim \pi_0(X/R)^\ff(R_n)$. Or on voit que
$\pi_0(X/R)^\ff$ n'a pas de sections sur $R$, si bien que la
composante formelle $(U_n)$ ne peut pas {\^e}tre approxim{\'e}e.
Les foncteurs $\pi_0(\sX/S)^\ff$ et $\Irr(\sX/S)^\ff$ ne sont
donc pas repr{\'e}sentables.



\no (3) Si $\sX$ est s{\'e}par{\'e} non propre, les espaces
tangents et des espaces d'obstructions de la th{\'e}orie de
d{\'e}formations de $\pi_0(\sX/S)^\ff$ ne sont pas de dimension
finie en g{\'e}n{\'e}ral. Par exemple, soit $X=X_0[\epsilon]$
comme dans l'exemple (1) et soit $Y=X_0\in\pi_0(\sX/S)^\ff(k)$.
Alors, l'espace tangent au foncteur de d{\'e}formations de
$Y$ est $\Hom_{\cO_X}(\cI_Y,\cO_X/\cI_Y)\simeq H^0(X_0,\cO_{X_0})$
et l'espace d'obstructions est
$\Ext^1_{\cO_X}(\cI_Y,\cO_X/\cI_Y)
\simeq \Ext^1_{\cO_{X_0}}(\cO_{X_0},\cO_{X_0})$
qui ne sont pas de dimension finie en g{\'e}n{\'e}ral.
\end{exems}

\begin{lemm} \label{lemm:pi_0_r_quasi_compact}
Pour tout $S$-champ alg{\'e}brique de pr{\'e}sentation finie
$\sX$, le foncteur $\pi_0(\sX/S)^\fr$
est localement de pr{\'e}sentation finie, quasi-compact et {\`a}
fibres finies.
\end{lemm}

\begin{demo}
Le fait que $\pi_0(\sX/S)^\fr$ est localement de pr{\'e}sentation
finie provient des r{\'e}sultats habituels de \cite{EGA}~IV.8.6
et de leur adaptation aux champs. Pour montrer qu'il est
quasi-compact et {\`a} fibres finies, on peut faire quelques
r\'eductions~: on peut supposer que $S$ est affine quitte \`a
localiser, qu'il est de type fini sur $\bZ$ puisque $\sX$ provient
par changement de base d'un tel sch{\'e}ma, et enfin que $S$ est
r{\'e}duit et irr{\'e}ductible, quitte {\`a} remplacer $S$ par une
de ses composantes irr{\'e}ductibles r{\'e}duites. Notons $\eta$ le
point g{\'e}n{\'e}rique de $S$. Par r{\'e}currence noeth{\'e}rienne,
il suffit de trouver un voisinage de $\eta$ au-dessus duquel
$\pi_0(\sX/S)^\fr$ est quasi-compact. On proc\`ede alors exactement
comme dans la d\'emonstration du lemme~\ref{lemme:rep_sur_un_ouvert}.
Les seules modifications sont le choix de l'extension finie
$K/k(\eta)$ qui doit \^etre telle que les sous-sch{\'e}mas ferm{\'e}s
r{\'e}duits des composantes connexes de $\sX_\eta\otimes K$ soient
g{\'e}om{\'e}triquement connexes et g\'eom\'etriquement r{\'e}duits.
Ensuite, on utilise \ref{prop:EGA_95}~(iv) et
\ref{theo:EGA_977}~(ii)-(iii) pour trouver un voisinage de $\eta$
au-dessus duquel on dipose de $r$ composantes connexes r{\'e}duites
relatives $\sC_i$ pour $\sX$. On trouve un morphisme surjectif
$S\times\{1,\dots,r\}\to\pi_0(\sX/S)^\fr$ ce qui
prouve notre assertion. Noter que contrairement \`a ce qu'il se
passe dans le lemme~\ref{lemme:rep_sur_un_ouvert}, on ne peut pas
prouver que ce morphisme est \'etale ou m\^eme plat.
\end{demo}



\subsection{Repr{\'e}sentabilit{\'e}}

\begin{theo} \label{theo:pi_0_cas_propre}
Soit $\sX$ un $S$-champ alg{\'e}brique propre et de pr{\'e}sentation
finie.
\begin{trivlist}
\itemn{1} Le foncteur $\pi_0(\sX/S)^\ff$ est repr{\'e}sentable
par un $S$-espace alg{\'e}brique formel localement de pr{\'e}sentation
finie et s{\'e}par{\'e}.
\itemn{2} Le foncteur $\pi_0(\sX/S)^\fr$
est repr{\'e}sentable par un $S$-sch{\'e}ma formel quasi-fini
et s{\'e}par{\'e}.
\end{trivlist}
\end{theo}


\begin{demo}
Les deux foncteurs sont localement de pr{\'e}sentation finie,
par les arguments habituels.
Il est clair, d'apr{\`e}s le crit{\`e}re valuatif, qu'ils
sont s{\'e}par{\'e}s sur $S$~: en effet, si $R$ est un anneau de
valuation discr{\`e}te de corps de
fractions~$K$, une \ccf\ $\sC_R$ de $\sX\times_S\Spec(R)$ est
uniquement d{\'e}termin{\'e}e comme l'adh{\'e}rence sch{\'e}matique
dans $\sX$ de sa fibre g{\'e}n{\'e}rique $\sC_K$.

Montrons que $\pi_0(\sX/S)^\ff$ est repr{\'e}sentable par un
$S$-espace alg{\'e}brique formel.
Soit $H_{\sX/S}$ le foncteur de Hilbert des
sous-champs ferm{\'e}s de $\sX$ sur~$S$~: c'est le foncteur $\Quot$
associ{\'e} au faisceau $\cO_\sX$, repr{\'e}sentable par un espace
alg{\'e}brique s{\'e}par{\'e}, localement de pr{\'e}sentation finie et
v{\'e}rifiant le crit{\`e}re valuatif de propret{\'e} (voir~\cite{Ol},
th.~1.5). Par ailleurs, pour un sous-champ ferm{\'e} $\sW\subset\sX$
de compl{\'e}mentaire ouvert $\sU$, il est {\'e}quivalent de dire que
les fibres $\sW_s$ ont un support ouvert dans $\sX_s$ ou que les fibres
de $\sU$ ont un support propre. D'apr{\`e}s le
lemme~\ref{lemme:lieu_de_proprete} ci-dessous, cette condition est
repr{\'e}sent{\'e}e par un sous-champ ouvert $K$ de $H_{\sX/S}$. Ensuite, on
note que le lieu des points de $S$ tels que la fibre $\sW_s$ est
g{\'e}om{\'e}triquement connexe est un ferm{\'e} de~$K$ (voir
\cite{EGA},~IV.15.5.9)~; notons $K_0$ le sous-espace ferm{\'e} r{\'e}duit que
ce ferm{\'e} d{\'e}termine. Un morphisme $f:T\to K$ d{\'e}finit une \ccf\
de $\sX_T$ sur~$T$ si et seulement si son image ensembliste tombe dans
$K_0$. Il est {\'e}quivalent de dire que l'image sch{\'e}matique $Z$
de $f$ v{\'e}rifie $Z_\red\subset K_0$, ou encore que l'id{\'e}al
$\cI$ de $K_0$ dans $K$ est inclus dans la racine de l'id{\'e}al
$\cJ$ de $Z$. Or $K$ est somme disjointe d'espaces
alg{\'e}briques de type fini $K^i$, de sorte que pour chaque $i$ on
obtient l'existence d'un entier $n_i$ tel que
$(\cI_{|K_i})^{n_i}\subset\cJ_{|K_i}$, c'est-{\`a}-dire qu'en restriction
{\`a} $K_i$, l'image de $f$ est dans le $n_i$-i{\`e}me
voisinage de $K_0\cap K^i$ dans $K^i$. Finalement
$\pi_0(\sX/S)^\ff$ est repr{\'e}sentable par la somme des
compl{\'e}t{\'e}s de $K^i$ le long de $K_0\cap K^i$.

Concernant le foncteur $\pi_0(\sX/S)^\fr$, on note d'abord que la
composante connexe universelle $\sC^u\to \pi_0(\sX/S)^\fr$
est un morphisme propre donc pur. On peut lui appliquer le
th{\'e}or{\`e}me~\ref{theo:CC_CI_pour_les_champs}(i) et
conclure que le foncteur $\pi_0(\sX/S)^\fr$ est repr{\'e}sentable
par le sous-espace ouvert de $\pi_0(\sX/S)^\ff$ au-dessus duquel
$\sC^u$ est {\`a} fibres g{\'e}om{\'e}triquement r{\'e}duites.
De plus $\pi_0(\sX/S)^\fr$ est quasi-compact d'apr{\`e}s le
lemme~\ref{lemm:pi_0_r_quasi_compact}, donc il est quasi-fini.
Il s'ensuit que cet espace alg{\'e}brique formel est en fait
un sch{\'e}ma formel~: {\'e}tant donn{\'e} que ses tronqu{\'e}s
sont quasi-finis et s{\'e}par{\'e}s sur $S$, ceci provient en
effet de~\cite{Kn}, chap.~II, cor.~6.17.
\end{demo}

Dans la preuve, on a utilis{\'e} le lemme suivant~:

\begin{lemm} \label{lemme:lieu_de_proprete}
Soit $\sX$ un $S$-champ alg{\'e}brique propre, de pr{\'e}sentation
finie, et $\sU$ un ouvert de $\sX$ de pr{\'e}sentation finie sur $S$.
Alors, le lieu $S_0$ des points $s\in S$ tels que $\sU_s$ est propre
sur $k(s)$ est un ouvert, et la restriction $\sU\times_S S_0\to S_0$
est propre.
\end{lemm}

\begin{demo}
Par les arguments standard, on peut supposer $S$ noeth{\'e}rien.
En suivant la lecture du d{\'e}but du paragraphe~15.7 de~\cite{EGA}~IV,
on voit que les r{\'e}sultats des num{\'e}ros~15.7.1 {\`a} 15.7.7 sont
valables avec pour seule modification de remplacer le crit{\`e}re
valuatif de propret{\'e} pour les sch{\'e}mas par son analogue pour les
champs (\cite{LMB} chapitre~7 et \cite{Ol}, introduction).

On peut donc se ramener au cas o{\`u} $S$ est le spectre d'un anneau
de valuation discr{\`e}te {\`a} corps r{\'e}siduel alg{\'e}briquement
clos (voir en particulier~\cite{EGA}~IV.15.7.5 et 15.7.7).
Il suffit de montrer que si la fibre $\sU_s$ au-dessus du point
ferm{\'e} $s\in S$ est propre sur $k(s)$ (on la suppose non vide,
sans quoi il n'y a rien {\`a} d{\'e}montrer), alors $\sU$ est
propre sur $S$. Puisque $S$ est strictement hens{\'e}lien, le
sch{\'e}ma $S'$ de la factorisation de Stein
$\sX\to S'=\St(\sX/S)\to S$ est somme disjointe d'un nombre
fini de $S$-sch{\'e}mas locaux finis $S'_1,\dots,S'_n$.
On peut remplacer $S$ par l'un des $S'_i$, puis par une composante
irr{\'e}ductible r{\'e}duite $T$ de $S'_i$ (\cite{EGA}~II.5.4.5)
surjective sur $S$ (sans quoi il n'y a rien {\`a} d{\'e}montrer)
et enfin par la normalisation $\tilde T$ puisque celle-ci est finie
sur $T$. On s'est donc ramen{\'e} au cas o{\`u} $S$ est le spectre d'un
anneau de valuation discr{\`e}te et $\sX\to S$ est {\`a} fibres
g{\'e}om{\'e}triquement connexes.

Par hypoth{\`e}se $\sU_s$ est propre sur $k(s)$, donc ouvert et
ferm{\'e} non vide dans $\sX_s$, donc $\sU_s=\sX_s$. Par suite,
le compl{\'e}mentaire $\sZ=\sX\setminus \sU$ est inclus dans la
fibre ouverte de $\sX$. Comme $\sZ$ est propre sur~$S$, ceci n'est
possible que si $\sZ=\emptyset$, ce qui montre que $\sU=\sX$.
\end{demo}

\begin{rema}
Contrairement au cas du foncteur des composantes
ouvertes avec $\sX$ \`a fibres g\'eom\'etriquement
r\'eduites, ici il n'existe pas en g\'en\'eral de morphisme
$\sX\to \pi_0(\sX/S)^\ff$ ou $\sX\to \pi_0(\sX/S)^\fr$.
Par exemple, soit $S$ le spectre d'un corps imparfait $k$
de caract\'eristique $p>0$, soit $t\not\in k^p$ et $C$
la courbe de type Fermat d'\'equation $x^p+y^p=tz^p$. On peut
voir que $\pi_0(C/k)^\fr\simeq\Spec(k(\sqrt[p]{t}))$ et il
n'y a pas de morphisme $C\to\pi_0(C/k)^\fr$.
\end{rema}

\begin{prop} \label{prop:comparaison_ouverte_fermee}
Soit $\sX$ un $S$-champ alg{\'e}brique de pr{\'e}sentation
finie, plat, {\`a} fibres g{\'e}om{\'e}\-triquement r{\'e}duites.
\begin{trivlist}
\itemn{1} On a
$\pi_0(\sX/S)^\ff=\pi_0(\sX/S)^\fr\subset\pi_0(\sX/S)$.
\itemn{2} On a un monomorphisme de foncteurs
$\Irr(\sX/S)^\ff\into\Irr(\sX/S)$.
\itemn{3} Si de plus $\sX$ est plat et pur sur $S$,
les monomorphismes
$$
\pi_0(\sX/S)^\ff\subset\pi_0(\sX/S)
\quad\mbox{et}\quad
\Irr(\sX/S)^\ff\into\Irr(\sX/S)
$$
sont des immersions ouvertes.
\end{trivlist}
\end{prop}

\begin{demo}
(1) Le fait que $\pi_0(\sX/S)^\ff=\pi_0(\sX/S)^\fr$ est clair.
Pour montrer que
$\pi_0(\sX/S)^\ff\subset\pi_0(\sX/S)$, il suffit d'observer
que toute \ccf\ $\sC\subset\sX$ est ouverte dans $\sX$.
On peut v{\'e}rifier ceci localement sur $S$, or comme le lieu
lisse de $\sC$ est $S$-dense, on peut supposer qu'on dispose
d'une section $g:S\to\sC$. Dans ce cas $\sC$ est {\'e}gale {\`a}
la composante connexe de $\sX$ le long de la section $g$, qui
est ouverte d'apr{\`e}s la proposition~\ref{comp_connexe_d_une_section}.

\medskip

\no (2) Soit $\sU$ le sous-champ ouvert de $\sX$ {\'e}gal au
lieu unicomposante. Si $\sI$ est une \cif\ de $\sX/S$,
le sous-champ $\sI\cap \sU$ est ouvert dans $\sX$. En effet,
en proc{\'e}dant comme dans (1) on voit qu'apr{\`e}s choix d'une
section locale $g:S\to\sI\cap \sU$, ce sous-champ est {\'e}gal
{\`a} la composante irr{\'e}ductible ouverte de $\sX$ le long
de la section $g$, qui est ouverte
(proposition~\ref{comp_irred_d_une_section}). Le morphisme
$\Irr(\sX/S)^\ff\into\Irr(\sX/S)$ est d{\'e}fini par
$\sI\mapsto \sI\cap\sU$. Comme de plus $\sI\cap\sU$ est plat
sur $S$ et dense dans $\sI$ fibre {\`a} fibre, il est
sch{\'e}matiquement dense dans $\sI$. Ainsi $\sI$ est {\'e}gal
{\`a} l'adh{\'e}rence sch{\'e}matique de $\sI\cap\sU$ dans $\sX$.
Ceci montre que le morphisme est un monomorphisme.

\medskip

\no (3) Soit $\sC^u$ la \cco\ universelle au-dessus de
$\pi_0(\sX/S)$. D'apr{\`e}s ce qui pr{\'e}c{\`e}de
$\pi_0(\sX/S)^\ff$ est repr{\'e}sen\-table par le
sous-espace de $\pi_0(\sX/S)$ des points o{\`u} l'adh{\'e}rence
sch{\'e}matique $\bar{\sC}{}^u$ de $\sC^u$ dans
$\sX\times \pi_0(\sX/S)$ est plate sur $S$ et {\`a} fibres
g{\'e}om{\'e}triquement connexes. Or, le lieu de platitude
de $\bar{\sC}{}^u\to S$ est un ouvert de $S$.
Restreignons-nous {\`a} cet ouvert et supposons donc
$\bar{\sC}{}^u$ plat sur $S$. On v{\'e}rifie
imm{\'e}diatement que $\bar{\sC}{}^u$ est pur sur $S$~:
en effet, la formation de l'adh{\'e}rence sch{\'e}matique
commute au changement de base plat de sorte qu'on peut
supposer que $S$ est local hens{\'e}lien~; comme $\sX$ est
pur sur $S$, pour tout point $x$ de $\bar{\sC}{}^u$ qui
est associ{\'e} dans sa fibre, l'adh{\'e}rence de $x$ (qui est
incluse dans $\bar{\sC}{}^u$ car celui-ci est ferm{\'e})
rencontre la fibre ferm{\'e}e de $\sX$ en un point
qui appartient donc {\`a} la fibre ferm{\'e}e de $\bar{\sC}{}^u$,
comme souhait{\'e}. Alors, le lieu des points o{\`u} la fibre
est g{\'e}om{\'e}triquement connexe est ouvert d'apr{\`e}s
le th{\'e}ror{\`e}me~\ref{theo:CC_CI_pour_les_champs}(ii),
ce qui conclut. On proc{\`e}de pareil pour montrer que
$\Irr(\sX/S)^\ff\into\Irr(\sX/S)$ est une immersion
ouverte.
\end{demo}

\begin{prop} \label{prop:propre_et_reduit}
Soit $\sX$ un $S$-champ alg{\'e}brique propre, plat, de
pr{\'e}sentation finie et {\`a} fibres g{\'e}om{\'e}triquement
r{\'e}duites.
\begin{trivlist}
\itemn{1} Soit $\sX\to \St(\sX/S)\to S$ la factorisation
de Stein. Alors, on a des isomorphismes
$$\St(\sX/S)\simeq\pi_0(\sX/S)^\ff=\pi_0(\sX/S)^\fr=\pi_0(\sX/S) \ .$$
\itemn{2} L'ouvert $\Irr(\sX/S)^\ff\subset\Irr(\sX/S)$
(voir~\ref{prop:comparaison_ouverte_fermee}(3)) est un sch{\'e}ma
{\'e}tale et s{\'e}par{\'e}.
\end{trivlist}
\end{prop}

\begin{demo}
(1) Notons qu'un champ propre est pur.
Le lemme~\ref{lemme:lieu_de_proprete} montre qu'une \cco\
$\sC\in\pi_0(\sX/S)$, qui est propre fibre {\`a} fibre, est propre
sur $S$. Ainsi, l'inclusion $\pi_0(\sX/S)^\ff\subset\pi_0(\sX/S)$
de~\ref{prop:comparaison_ouverte_fermee}(3) est
en fait une {\'e}galit{\'e}.
Notons maintenant $S':=\St(\sX/S)$. Pour tout point $s\in S$,
l'alg{\`e}bre $H^0(\sX_s,\cO_{\sX_s})$ est {\'e}tale car la fibre
$\sX_s$ est g{\'e}om{\'e}triquement r{\'e}duite. Il en d{\'e}coule
que $f:\sX\to S$ est cohomologiquement plat en dimension $0$ et
que $f_*\cO_{\sX}$ est une $\cO_S$-alg{\`e}bre {\'e}tale
(\cite{EGA}~III.7.8.6 et 7.8.7). Ainsi $S'\to S$
est {\'e}tale, donc le morphisme $\sX\to S'$ est plat. Soit $s'$
un point g{\'e}om{\'e}trique de $S'$ et $s$ son image
dans $S$. D'apr{\`e}s les propri{\'e}t{\'e}s de la factorisation
de Stein, les fibres de $\sX\to S'$ sont des composantes connexes
de $\sX_s$ de sorte que finalement $\sX$ est une \ccf\ de
$\sX\times_S S'/S'$ correspondant {\`a} un morphisme
$S'\to\pi_0(\sX/S)^\ff$. Ce morphisme est un morphisme entre
deux sch{\'e}mas finis {\'e}tales, qui est un isomorphisme fibre {\`a}
fibre (noter que la formation de $\St(\sX/S)$ commute au
changement de base), donc c'est un isomorphisme.

\medskip

\no (2) Le fait que $\Irr(\sX/S)^\ff$ est s\'epar\'e provient du
fait que sur un anneau de valuation discr{\`e}te de base, une
\cif\ est uniquement d\'etermin\'ee comme adh{\'e}rence
sch{\'e}matique dans $\sX$ de sa fibre g{\'e}n{\'e}rique. Comme
$\Irr(\sX/S)^\ff$ est \'etale et s\'epar\'e, c'est alors un sch\'ema.
\end{demo}

La fonctorialit\'e en $\sX$ pour $\pi_0(\sX/S)^\ff$ est
\'evidemment moins bonne que pour $\pi_0(\sX/S)$. Nous
terminons cette section avec deux cas simples de covariance
et de contravariance~:

\begin{prop}
Soient $\sX$ et $\sY$ des $S$-champs alg{\'e}briques plats et
de pr{\'e}sentation finie.
\begin{trivlist}
\itemn{1} Tout morphisme fini \'etale $f:\sX\to\sY$ induit
un morphisme $f_*:\pi_0(\sX/S)^\ff\to\pi_0(\sY/S)^\ff$.
\itemn{2} Tout morphisme fppf \`a fibres g\'eom\'etriquement
connexes $f:\sX\to\sY$ induit un morphisme
$f^*:\pi_0(\sY/S)^\ff\to\pi_0(\sX/S)^\ff$.
\end{trivlist}
\end{prop}

\begin{demo}
Soit $f:\sX\to\sY$ un morphisme avec $\sX,\sY$ plats et de
pr{\'e}sentation finie.

\no (1) Le morphisme covariant $f_*$ envoie une \ccf\ $\sC$ sur
l'image sch\'ematique $\sD:=f(\sC)$. Comme $f$
est \'etale, $\sD$ est plat sur $S$. Comme
$g:\sC\into\sX\to\sY$ est non ramifi\'e, le $\cO_\sY$-module
$g_*\cO_\sC$ est localement monog\`ene, donc son annulateur
est \'egal \`a $\Fitt_0(g_*\cO_\sC)$, le $0$-i\`eme id\'eal
de Fitting (\cite{Ei}, prop.~20.7), et ceci reste vrai
apr\`es tout changement de base. Comme la formation des
id\'eaux de Fitting commute au changement de base, il en va
de m\^eme pour la formation de $\sD$. Enfin, comme $f$
est ouvert et ferm\'e, le sous-champ $\sD\subset \sY$ est
\`a fibres g\'eom\'etriquement connexes, ouvertes et
ferm\'ees donc finalement c'est bien une \ccf\ de $\sY/S$.

\no (2) Le morphisme contravariant $f^*$ envoie une \ccf\
$\sD\subset\sY$ sur $\sC:=f^{-1}(\sD)$. Il est clair que $\sC$
est plat et de pr\'esentation finie sur $S$, et comme $f$
est universellement submersif \`a fibres g\'eom\'etriquement
connexes, ses fibres g\'eom\'etriques
qui sont ouvertes et ferm\'ees sont aussi connexes.
\end{demo}

\subsection{Un contre-exemple} \label{ss_section:contre_ex}

Dans cette sous-section, nous montrons que pour un $S$-champ
alg{\'e}brique propre, plat et de pr{\'e}sentation finie,
le foncteur $\Irr(\sX/S)^\ff$ n'est pas en g{\'e}n{\'e}ral
repr{\'e}sentable par un $S$-espace alg{\'e}brique formel.
Le contre-exemple est donn\'e dans~\ref{adoc:contre_ex_coniques}.
Il est bas{\'e} sur la
propri{\'e}t{\'e}~\ref{prop:lieu_vs_foncteur} ci-dessous,
qui est peut-{\^e}tre bien connue mais dont je n'ai trouv{\'e}
mention nulle part dans la litt{\'e}rature. Notons que les
crit\`eres d'Artin montrent assez facilement que le foncteur
$\Irr(\sX/S)^\ff$ de l'exemple~\ref{adoc:contre_ex_coniques}
n'est pas repr{\'e}sentable par un $S$-espace
alg{\'e}brique. Malheureusement, les crit\`eres d'Artin ne
s'\'etendent pas aux espaces alg\'ebriques formels (notamment,
l'axiome d'approximation des \'el\'ements formels fait d\'efaut).
C'est pourquoi nous avons recours \`a une id\'ee diff\'erente.

Soit $P=P(\sX/k)$ une propri{\'e}t{\'e} des champs alg{\'e}briques $\sX$
de type fini sur un corps $k$ qui est invariante par extension
du corps de base, au sens o{\`u} $\sX$ v{\'e}rifie $P$ si et seulement
si $\sX\otimes_k \ell$ v{\'e}rifie $P$, pour toute extension de
corps $\ell/k$. On dira aussi que $P$ est une propri{\'e}t{\'e}
{\em g{\'e}om{\'e}trique}. {\`A} tout champ alg{\'e}brique $\sX$ de
pr{\'e}sentation finie sur un sch{\'e}ma $S$, on peut associer~:
\begin{trivlist}
\itemn{1} l'ensemble $E_P=E_P(\sX/S)$ des $s\in S$ tels que la
fibre de $\sX$ en $s$ v{\'e}rifie $P$, et
\itemn{2} le foncteur $F_P=F_P(\sX/S)$ sur la cat{\'e}gorie des
$S$-sch{\'e}mas d{\'e}fini par
$$
F_P(T)=
\left\{
\begin{array}{cl}
\{\emptyset\} &
\mbox{si les fibres de $X\times_S T\to T$ v{\'e}rifient $P$,} \\
\emptyset & \mbox{sinon.} \\
\end{array}
\right.
$$
\end{trivlist}
On appellera $E_P$ le {\em lieu indicateur} et $F_P$ le
{\em foncteur indicateur} de la propri{\'e}t{\'e} $P$ pour $\sX/S$.

\begin{prop} \label{prop:lieu_vs_foncteur}
Avec les notations ci-dessus~:
\begin{trivlist}
\itemn{1} $F_P$ est repr{\'e}sentable par un sch{\'e}ma si et
seulement si $E_P$ est ouvert dans $S$.
\itemn{2} $F_P$ est repr{\'e}sentable par un sch{\'e}ma formel si et
seulement si $E_P$ est localement ferm{\'e} dans~$S$.
\end{trivlist}
\end{prop}

\begin{remas}
(1) Si $S$ est un espace alg{\'e}brique au lieu d'un
sch{\'e}ma, le m{\^e}me {\'e}nonc{\'e} est valable (avec la m{\^e}me
preuve) en rempla{\c c}ant <<~sch{\'e}ma~>> par
<<~espace alg{\'e}brique~>> et <<~sch{\'e}ma formel~>> par
<<~espace alg{\'e}brique formel~>>.

\no (2) La preuve utilise trois petits lemmes
(\ref{lemm:mono_vers_un_trait},
\ref{lemm:specialisations_successives} et
\ref{lemm:non_existence_de_mono}) qui seront {\'e}tablis
ci-dessous.
\end{remas}

\begin{demo}
La formation de $E_P$ et de $F_P$ commute aux changements de
base $S'\to S$. Comme les termes des {\'e}quivalences {\`a}
d{\'e}montrer sont des propositions de nature locale sur $S$,
on peut supposer $S$ affine. Par ailleurs, comme $\sX$ est de
pr{\'e}sentation finie sur $S$, il en va de m{\^e}me pour le
foncteur $F_P$, et on peut donc supposer que $S$ est le spectre
d'un anneau noeth{\'e}rien.

\medskip

\no (1) Si $E_P$ est ouvert dans $S$, alors c'est un sous-sch{\'e}ma
de $S$ et il est clair qu'il repr{\'e}sente $F_P$. R{\'e}ciproquement,
si $F_P$ est repr{\'e}sentable par un sch{\'e}ma, alors la partie $E_P$,
qui est l'image de $F_P\to S$, est constructible. Pour montrer
que $E_P$ est ouvert, il suffit de montrer qu'il est stable par
g{\'e}n{\'e}risation. Pour cela, on peut remplacer $S$ par le spectre
d'un anneau de valuation discr{\`e}te $(R,m)$, dont le point ferm{\'e}
est dans $E_P$, et il faut montrer que $E_P=S$. Or par hypoth{\`e}se
$F_P(R/m^n)=\{\emptyset\}$ pour tout $n\ge 1$. Comme aucune immersion
ferm{\'e}e $f:X\to S$ distincte de l'identit\'e ne se factorise {\`a}
travers tous les voisinages
infinit{\'e}simaux $\Spec(R/m^n)$, le
lemme~\ref{lemm:mono_vers_un_trait} ci-dessous montre que l'image
de $F_P\to S$ est n{\'e}cessairement {\'e}gale {\`a} $S$, donc $E_P=S$.

\medskip

\no (2) Si $E_P$ est localement ferm{\'e} dans $S$, il est ferm{\'e}
dans un ouvert $U\subset S$. Pour montrer que $F_P$ est
repr{\'e}sentable par un sch{\'e}ma formel, on peut remplacer $S$ par
$U$ et supposer que $E_P$ est ferm{\'e}. On note encore $E_P$
le sous-sch{\'e}ma ferm{\'e} r{\'e}duit de $S$ de support $E_P$.
Un morphisme $f:T\to S$ d{\'e}finit un point de $F_P$ si et
seulement si son image ensembliste tombe dans $E_P$. Il est
{\'e}quivalent de dire que l'image sch{\'e}matique $Z$ de $f$
v{\'e}rifie $Z_\red\subset E_P$, ou encore que l'id{\'e}al
$\cI$ de $E_P$ dans $K$ est inclus dans la racine de l'id{\'e}al
$\cJ$ de $Z$. Or $E_P$ est noeth{\'e}rien, donc il existe un entier
$n$ tel que $\cI^n\subset\cJ$, ce qui signifie que~$f$ se factorise
par le $n$-i{\`e}me voisinage de $E_P$ dans $S$. Finalement
$F_P$ est repr{\'e}sentable par le compl{\'e}t{\'e} de~$S$ le long
de $E_P$.

R{\'e}ciproquement, supposons que $F_P$ est repr{\'e}sentable
par un sch{\'e}ma formel. Il suffit de montrer que $E_P$ est
ouvert dans son adh{\'e}rence. On peut donc remplacer $S$ par
l'adh{\'e}rence r{\'e}duite de $E_P$ dans $S$ et~$X$ par sa
restriction {\`a} cette adh{\'e}rence. Alors $E_P$ est dense dans
$S$ et on doit montrer qu'il est ouvert. Il suffit de montrer
qu'il est stable par g{\'e}n{\'e}risation. Soit $x_2\in E_P$ et
$x_1\in S$ une g{\'e}n{\'e}risation de $x$. Comme $E_P$ est dense
dans $S$, il existe une g{\'e}n{\'e}risation $x_0$ de $x_1$ qui
appartient {\`a} $E_P$. D'apr{\`e}s le
lemme~\ref{lemm:specialisations_successives}, il existe un
sch{\'e}ma $\Sigma$, spectre d'un anneau de valuation de rang $2$,
de points $\sigma_0\leadsto \sigma_1\leadsto \sigma_2$ (o{\`u}
$\sigma_i$ est le point de codimension $i$ de $\Sigma$) et
un morphisme $f:\Sigma\to S$ tel que $f(s_i)=x_i$ pour
$i=0,1,2$. En faisant le changement de base $\Sigma\to S$, on
se ram{\`e}ne au cas o{\`u}~$S$ est le spectre d'un anneau de
valuation de rang~$2$. Notons $F_0$ la fibre sp{\'e}ciale de $F_P$,
qui est un $S$-sch{\'e}ma de pr{\'e}sentation finie. Supposons
maintenant que $x_1\not\in E_P$, alors l'image du monomorphisme
$F_0\to F_P\to S$ est $\{x_0,x_2\}$. Si l'on choisit un trait
$T$ et un morphisme $T\to S$ d'image $\{x_0,x_2\}$, on voit
que la restriction de $F_P$ {\`a} $T$ est repr{\'e}sentable par $T$
lui-m{\^e}me et est donc connexe. Ceci montre que $F_P$ et $F_0$
sont connexes. D'apr{\`e}s le lemme~\ref{lemm:non_existence_de_mono},
ceci est exclu. Il s'ensuit que $x_1\in E_P$ donc $E_P$ est
ouvert, ce qui conclut la preuve de la proposition.
\end{demo}

\begin{lemm} \label{lemm:mono_vers_un_trait}
Soit $f:X\to S$ un monomorphisme de sch{\'e}mas tel que $S$
est un trait et $f(X)$ est le point ferm{\'e}. Alors $f$ est
une immersion ferm{\'e}e.
\end{lemm}

\begin{demo}
Comme $f$ est un monomorphisme, $X$ est r{\'e}duit {\`a} un point
et il est donc affine, d'anneau $A$ artinien. La restriction
de $f$ au-dessus du point ferm{\'e} est un monomorphisme d'un
sch{\'e}ma non vide {\`a} valeurs dans le spectre d'un corps donc
c'est un isomorphisme. D'apr{\`e}s le lemme de Nakayama, le
morphisme $\Gamma(S,\cO_S)\to A$ est donc surjectif.
\end{demo}

\begin{lemm} \label{lemm:specialisations_successives}
Soient $s_0\leadsto s_1\leadsto\dots\leadsto s_n$ des points d'un
sch{\'e}ma localement noeth{\'e}rien $S$ tels que $s_i$ est une
sp{\'e}cialisation de $s_{i-1}$ pour tout $i=1,\dots,n$. Alors il
existe un sch{\'e}ma $T$, spectre d'un anneau de valuation de rang
$n$, de points not{\'e}s $t_0\leadsto t_1\leadsto\dots\leadsto t_n$
o{\`u} $t_i$ est l'unique point de codimension~$i$, et un morphisme
$f:T\to S$ tel que $f(t_i)=s_i$ pour tout $i$.
\end{lemm}

\begin{demo}
On peut remplacer $S$ par l'adh{\'e}rence de $s_0$ puis par son
sous-sch{\'e}ma r{\'e}duit et donc supposer $S$ int{\`e}gre de corps de
fonctions {\'e}gal au corps r{\'e}siduel de $s_0$. On peut ensuite
remplacer $S$ par un ouvert affine contenant $s_n$ et on se ram{\`e}ne
ainsi au cas o{\`u} $S$ est le spectre d'un anneau int{\`e}gre noeth{\'e}rien
$A$, de corps de fractions $K=k(S)$. Notons $p_i\subset A$ l'id{\'e}al
premier correspondant au point $s_i$ et $A_{p_i}$ son anneau local,
avec $p_0=(0)\subset p_1\subset\dots\subset p_n\subset A$ et
$A_{p_n}\subset\dots\subset A_{p_1}\subset A_{p_0}=K$. On consid{\`e}re~:
\begin{trivlist}
\itemn{1}
une valuation $v_1:K\to \bZ$ dont l'anneau de valuation domine
$A_{p_1}$. Alors $A_{p_1}/p_1$ s'identifie {\`a} un sous-anneau
du corps r{\'e}siduel $k(v_1)$, de corps de fractions $k(v_1)$~;
\itemn{2}
une valuation $v_2:k(v_1)\to \bZ$ dont l'anneau de valuation domine
$A_{p_2}/p_1$. Alors $A_{p_2}/p_2$ s'identifie {\`a} un sous-anneau
du corps r{\'e}siduel $k(v_2)$, de corps de fractions $k(v_2)$~;
\end{trivlist}
et ceci jusqu'{\`a}~:
\begin{trivlist}
\itemn{$n$}
une valuation $v_n:k(v_{n-1})\to\bZ$ dont l'anneau de valuation domine
$A_{p_n}/p_{n-1}$.
\end{trivlist}
On consid{\`e}re alors la valuation lexicographique
$v=(v_1,\dots,v_2):K\to\bZ^n$ associ{\'e}e aux $v_i$ et {\`a} un choix
d'uniformisantes $\pi_i\in k(v_{i-1})$ telles que $v_i(\pi_i)=1$.
Elle est d{\'e}finie pr{\'e}cis{\'e}ment ainsi~: si $x\in K\setminus\{0\}$,
on note $i_1:=v_1(x)$ et $x_1$ la classe r{\'e}siduelle de
$\pi_1^{-i_1}x$ dans $k(v_1)$, puis $i_2:=v_2(x_1)$ et $x_2$ la classe
r{\'e}siduelle de $\pi_2^{-i_2}x_1$ dans $k(v_2)$, etc. On pose alors
$v(x)=(i_1,\dots,i_n)$. On note enfin
$V=\{x\in K,\,v(x)\ge (0,\dots,0)\}$ l'anneau de valuation de $v$
et $T=\Spec(V)$. Par construction, on a un morphisme $f:T\to S$
qui satisfait aux conditions de l'{\'e}nonc{\'e}.
\end{demo}

\begin{lemm} \label{lemm:non_existence_de_mono}
Soit $S$ le spectre d'un anneau de valuation de rang $2$ et
$s_0\leadsto s_1\leadsto s_2$ ses points. Alors, il n'existe
pas de monomorphisme de pr{\'e}sentation finie $f:X\to S$ tel que
$X$ est un sch{\'e}ma connexe et $f(X)=\{s_0,s_2\}$.
\end{lemm}

\begin{demo}
Soit $i$ l'immersion ferm{\'e}e $X_\red\into X$, quitte {\`a}
consid{\'e}rer $f\circ i$ {\`a} la place de $f$ on peut supposer que
$X$ est r{\'e}duit. Comme $X$ est connexe, la pr{\'e}image par $f$ du
point ferm{\'e} $s_2$ est un point $y$ ferm{\'e} dans $X$ mais non
ouvert. Il s'ensuit que tout ouvert affine de $X$ contenant $y$
est {\'e}gal {\`a} $X$, donc $X$ est affine. La topologie de $X$ est
celle d'un trait, en particulier son anneau $B$ est int{\`e}gre.
Soit $V$ l'anneau de valuation dont $S$ est le spectre et $K$ son
corps de fractions. La restriction de $X\to S$ au-dessus du
point g{\'e}n{\'e}rique est un monomorphisme d'un sch{\'e}ma non vide {\`a}
valeurs dans $\Spec(K)$ donc c'est un isomorphisme, donc le corps
de fractions de $B$ est $K$. De plus $B$ est sans torsion comme
$V$-module, donc plat sur $V$. Comme $X\to S$ est de pr{\'e}sentation
finie et plat il est ouvert, contradiction.
\end{demo}

\begin{adoc}{Le contre-exemple} \label{adoc:contre_ex_coniques}
L'exemple suivant m'a {\'e}t{\'e} sugg{\'e}r{\'e} par Angelo Vistoli.
Il montre que $\Irr(X/S)^\ff$ ne poss{\`e}de pas d'aussi bonnes
propri{\'e}t{\'e}s de repr{\'e}sentabilit{\'e} que $\pi_0(\sX/S)^\ff$.
Sur un corps $k$ de caract{\'e}ristique diff{\'e}rente de $2$,
nous consid{\'e}rons l'espace modulaire des coniques planes
$S:=\bP^5=\Proj(k[a,b,c,d,e,f])$ et la conique universelle
$X\subset\bP^2\times\bP^5$ d'{\'e}quation $q(x,y,z)=0$ o{\`u}
$q(x,y,z)=ax^2+by^2+cz^2+dxy+exz+fyz$. Il y a trois types de
coniques~:

\medskip

\begin{itemize}
\item les coniques lisses correspondent {\`a} l'ouvert
$U=\{\disc(q)\ne 0\}$~;
\item les droites doubles vivent dans le
ferm{\'e} $F$ image du morphisme $(\bP^2)^\vee\simeq\bP^1\to S$
qui envoie une droite d'{\'e}quation $\ell=0$ sur la conique
d'{\'e}quation $\ell^2=0$~;
\item les coniques singuli{\`e}res r{\'e}ductibles forment la
partie localement ferm{\'e}e $(X\setminus U)\setminus F$.
\end{itemize}

\medskip

On note que l'ensemble des $s\in S$ tels que $X_s$ est
g{\'e}om{\'e}triquement irr{\'e}ductible n'est pas localement ferm{\'e}.
Posons $F:=\Irr(X/S)^\ff$ et montrons que $F$ n'est pas
repr{\'e}sentable par un $S$-espace alg{\'e}brique formel.
Soit $H\to S$ le sch{\'e}ma de Hilbert
des sous-sch{\'e}mas ferm{\'e}s de~$X$ et soit $H_0\subset H$
l'ouvert et ferm{\'e} contenant le sous-sch{\'e}ma ferm{\'e}
maximal $Z=X$. C'est la composante du sch{\'e}ma de Hilbert
indic{\'e}e par le polyn{\^o}me de Hilbert maximal $P(n)=2n+1$,
et on a $H_0\simeq S$. Soit $F_0\subset F$ l'ouvert et
ferm{\'e} pr{\'e}image. Il est clair que $F_0$ est le sous-foncteur
de $H_0$ indicateur du lieu o{\`u} les fibres de $X\to S$ sont
g{\'e}om{\'e}triquement irr{\'e}ductibles. Comme l'ensemble des points
$s\in S$ tels que $Z_s$ est g{\'e}om{\'e}triquement irr{\'e}ductible
n'est pas localement ferm{\'e}, il d{\'e}coule de la
proposition~\ref{prop:lieu_vs_foncteur} que $F_0$ n'est pas
repr{\'e}sentable par un $S$-espace alg{\'e}brique formel.
A fortiori, $F$ n'est pas repr{\'e}sentable par un $S$-espace
alg{\'e}brique formel. \qed
\end{adoc}

\subsection{Exemple~: modules des courbes admettant une action}

Soit $G$ un groupe fini, $\gamma$ son cardinal et
$S=\Spec(\bZ[1/\gamma])$ qui sert de sch\'ema de base. On fixe
un entier $g\ge 2$ et on consid\`ere le champ $\sM_g$ des
courbes de genre $g$. C'est un champ de Deligne-Mumford lisse
et de dimension $3g-3$.

\begin{prop}
Soit $\sM_g(G)$ le sous-$S$-champ de $\sM_g$ des courbes qui admettent
une action fid{\`e}le de $G$.
\begin{trivlist}
\itemn{1} Le champ $\sM_g(G)$ est un sous-champ ferm{\'e}, que l'on
munit de la structure de sous-champ alg{\'e}brique r{\'e}duit.
Il est plat, de pr\'esentation finie sur $\bZ[1/\gamma]$, \`a
fibres g\'eom\'etriquement r\'eduites.
\itemn{2} La normalisation $\tilde\sM_g(G)$ de $\sM_g(G)$ est un champ
alg{\'e}brique lisse sur $\bZ[1/30\gamma]$.
\end{trivlist}
\end{prop}

La preuve utilise le lemme~4.1 de \cite{MSSV} qui montre plus
pr\'ecis\'ement qu'en fait, en dehors d'une liste explicite de
$10$ groupes, le r\'esultat de (2) vaut aussi sur $\bZ[1/2\gamma]$.

\begin{demo}
(1) Notons $\sH$ le champ de Hurwitz classifiant les paires
$(C,\phi)$ o{\`u} $C$ est une courbe de genre $g$ et
$\phi:G\hookrightarrow \Aut(C)$ est un monomorphisme de sch\'emas
en groupes. C'est un champ alg{\'e}brique de Deligne-Mumford lisse
sur ${\mathbb Z}[1/\gamma]$, non {\'e}quidimensionnel~: les dimensions
des diff{\'e}rentes composantes connexes d{\'e}pendent de la
ramification de l'action du groupe $G$. Le morphisme $f:\sH\to\sM_g$
donn\'e par l'oubli de l'action est fini (repr\'esentable) et
non ramifi\'e.
Ainsi $\sM_g(G)$, qui est l'image de $f$, ou encore son image
sch\'ematique, est ferm\'e. Les composantes irr\'eductibles de
$\sM_g(G)$ sont images de composantes irr\'eductibles de $\sH$,
et en particulier dominent $\Spec(\bZ[1/\gamma])$. Ceci montre
que $\sM_g(G)$ est plat sur $\bZ[1/\gamma]$. Par ailleurs, comme
$f$ est non ramifi\'e, le $\cO_{\sM_g}$-module $f_*\cO_\sH$ est
localement engendr\'e par un \'el\'ement, ce qui montre que
son annulateur est \'egal \`a $\Fitt_0(f_*\cO_\sH)$, le
$0$-i\`eme id\'eal de Fitting (\cite{Ei}, prop.~20.7). De plus,
ceci reste vrai apr\`es tout changement de base. Comme
la formation des id\'eaux de Fitting commute au changement de
base, il en va de m\^eme pour l'image (sch\'ematique) de $f$.
Alors, comme les fibres de $\sH$ sur $S$ sont g\'eom\'etriquement
r\'eduites, la m\^eme chose est vraie pour les fibres de $\sM_g(G)$.

\medskip

\no (2) On peut remplacer $\sM_g(G)$ par une de ses composantes
irr{\'e}ductibles. Notons alors $\eta$ le point g{\'e}n{\'e}rique
et $G'$ le groupe d'automorphisme de la courbe correspondante~;
on a donc $G\subset G'$ et il est clair que $\sM_g(G)=\sM_g(G')$.
De plus on sait d'apr\`es \cite{MSSV}, lemma~4.1 que les premiers
qui divisent $[G':G]$ sont dans $\{2,3,5\}$ (dans {\em loc. cit.},
ce r\'esultat est \'enonc\'e en caract\'eristique $0$, mais la
lecture de la preuve montre que tout est valable en
caract\'eristique premi\`ere \`a l'ordre de $G$). Ainsi, quitte
{\`a} se restreindre \`a $\bZ[1/30\gamma]$ et \`a remplacer $G$
par $G'$, on peut supposer que le groupe d'automorphismes de la
courbe g{\'e}n{\'e}rique est exactement $G$. Dans la suite, nous
supposons cette condition r\'ealis\'ee.

Le champ $\sH$ est muni d'une action {\`a} gauche du groupe
$\Aut(G)$ par <<~torsion des actions~>> d{\'e}finie de la
mani{\`e}re suivante~: si $\alpha$ est un automorphisme de $G$,
on pose $\alpha.(C,\phi)=(C,\phi\circ\alpha^{-1})$.
Le morphisme $f:\sH\to\sM_g$ est clairement invariant sous $\Aut(G)$,
de sorte qu'il induit un morphisme fini surjectif
$f':\sH/\Aut(G)\to\sM_g(G)$. Comme le groupe d'automorphismes
de la courbe g{\'e}n{\'e}rique est $G$, ce morphisme est
birationnel. Comme de plus le morphisme de quotient
$\sH\to \sH/\Aut(G)$ est {\'e}tale, le champ $\sH/\Aut(G)$ est
lisse sur $S$ donc normal.
D'apr{\`e}s le th{\'e}or{\`e}me principal de Zariski, $f'$
s'identifie {\`a} la normalisation de $\sM_g(G)$.
\end{demo}

Dans la terminologie en vigueur, ce r\'esultat montre que
$\sM_g(G)$ est \'equinormalisable, et m\^eme
\'equid\'esingularisable~:

\begin{coro} \label{normalisation_et_changement_de_base_normal}
Soit $S=\Spec(\bZ[1/30\gamma])$ et $S'\to S$ un changement de base avec
$S'$ normal. Alors $\tilde\sM_g(G)\times_S S'$ est la normalisation de
$\sM_g(G)\times_S S'$. Les fibres de $\tilde\sM_g(G)\to S$ sont
les normalisations des fibres de $\sM_g(G)$.
\end{coro}

\begin{demo}
Notons d'abord que $\sM_g(G)\times_S S'$ (resp. $\tilde\sM_g(G)\times_S S'$)
est {\`a} fibres g{\'e}om{\'e}triquement r{\'e}duites (resp. est lisse) sur
$S'$ normal, donc il est r{\'e}duit (resp. normal).
Le sous-champ ouvert $\sU\subset \sM_g(G)$ au-dessus duquel le
morphisme $f':\sH/\Aut(G)\to\sM_g(G)$ est plat sur $S$ et dense
fibre {\`a} fibre, donc $S$-universellement dense
(\cite{EGA}~IV.11.10.9). Il s'ensuit que le morphisme
$\tilde\sM_g(G)\times_S S'\to\sM_g(G)\times_S S'$
est un isomorphisme au-dessus de $\sU\times_S S'$, donc
 birationnel. Ce morphisme est aussi quasi-fini,
s{\'e}par{\'e}, surjectif, donc c'est la normalisation de
$\sM_g(G)\times_S S'$.
\end{demo}

\begin{coro} \label{coro:Irr_de_MgG_fini_etale}
Soit $S=\Spec(\bZ[1/30\gamma])$. Le foncteur $\Irr(\sM_g(G)/S)$
est repr{\'e}sentable par un $S$-sch{\'e}ma fini {\'e}tale.
\end{coro}

\begin{demo}
Puisque la normalisation $\pi:\tilde\sM_g(G)\to\sM_g(G)$ est
$S$-birationnelle, elle induit un isomorphisme
$\Irr(\tilde\sM_g(G)/S)\simeq\Irr(\sM_g(G)/S)$
(corollaire~\ref{coro:covariance_en_X}). Comme de plus
$\tilde\sM_g(G)/S$ est normal sur $S$, on a bien s{\^u}r
$\Irr(\tilde\sM_g(G)/S)=\pi_0(\tilde\sM_g(G)/S)$. Pour finir,
il est connu que le champ $\tilde\sM_g(G)=\sH/\Aut(G)$ admet
une compactification lisse $\bar\sH$ dans laquelle il est
$S$-dense (voir \cite{BR}, section~6.3). Il en d{\'e}coule
que $\pi_0(\tilde\sM_g(G)/S)$ est isomorphe {\`a}
$\pi_0(\bar\sH/S)$ et ce dernier foncteur est
repr{\'e}sentable par un sch{\'e}ma fini {\'e}tale,
d'apr{\`e}s la proposition~\ref{prop:propre_et_reduit}.
\end{demo}

\appendix

\section{Propri{\'e}t{\'e}s constructibles pour les champs alg{\'e}briques}
\label{Annexe_constructiblite}

Dans cette annexe, nous rappelons quelques r{\'e}sultats
de \cite{EGA}~IV concernant la constructiblit{\'e} de
certaines parties dans des champs alg{\'e}briques. Les
d{\'e}monstrations, {\'e}crites dans le cas des sch{\'e}mas,
s'adaptent de mani{\`e}re {\`a} peu pr{\`e}s imm{\'e}diate au cas
des morphismes de champs alg{\'e}briques, pourvu que l'on
soit suffisamment soigneux dans les {\'e}nonc{\'e}s. En g{\'e}n{\'e}ral,
les adaptations n{\'e}cessaires reviennent {\`a} remplacer,
lorsque n{\'e}cessaire, l'utilisation de recouvrements par
des sch{\'e}mas affines ouverts {\em de Zariski} par des
recouvrements par des sch{\'e}mas affines ouverts pour la
topologie {\em lisse}~; ou alors {\`a} appliquer les r{\'e}sultats
de \cite{EGA} pour des groupo{\"\i}des $X_1\rightrightarrows X_0$,
c'est-{\`a}-dire des paires de morphismes satisfaisant certaines
conditions, au lieu de les appliquer simplement pour des
sch{\'e}mas~; ou {\`a} utiliser d'autres techniques du m{\^e}me
genre, classiques lorsqu'on manipule des champs alg{\'e}briques.

Ci-dessous, nous r{\'e}unissons un ensemble de r{\'e}sultats qui
incluent d'une part ceux dont nous avons besoin dans le corps
de l'article, et d'autre part ceux qui sont n{\'e}cessaires {\`a}
la preuve (adapt{\'e}e de \cite{EGA}~IV) des pr{\'e}c{\'e}dents.
Apr{\`e}s chaque {\'e}nonc{\'e}, nous indiquons chaque fois que cela
est utile les modifications {\`a} apporter {\`a} la preuve de
\cite{EGA} pour passer des sch{\'e}mas aux champs alg{\'e}briques.

\begin{appadoc}{Fibres des morphismes de champs alg{\'e}briques}
Soit $f:\sX\to\sS$ un morphisme de champs alg{\'e}briques et
$s\in |\sS|$ un point. M{\^e}me lorsque $s$ poss{\`e}de un
{\em corps r{\'e}siduel $k(s)$} bien d{\'e}fini, il n'existe
pas en g{\'e}n{\'e}ral de morphisme $\Spec(k(s))\to\sS$, de sorte
que la notion de {\em fibre de $f$ au point $s$} n'est pas
aussi imm{\'e}diate que dans le cas des morphismes de sch{\'e}mas.
Pour prendre ce fait en compte, nous utiliserons la convention
de terminologie suivante.

Consid{\'e}rons une propri{\'e}t{\'e} de la forme
$P=P(\sX/k,\sF,Z)$ portant sur des donn{\'e}es compos{\'e}es d'un
champ alg{\'e}brique $\sX$ sur un corps $k$, un
$\cO_\sX$-module $\sF$ et une partie $Z\subset |\sX|$.
(De la m{\^e}me fa{\c c}on, les consid{\'e}rations qui suivent
sont valables pour des propri{\'e}t{\'e}s mettant en jeu
un nombre fini de champs alg{\'e}briques, de modules ou de
parties sur ces champs, de morphismes entre ces champs...
ou une partie seulement de ces donn{\'e}es.) On s'int{\'e}ressera
principalement
{\`a} des propri{\'e}t{\'e}s ind{\'e}pendantes du corps de base
au sens o{\`u}, pour toute extension de corps $\ell/k$,
$P(\sX/k,\sF,Z)$ est vraie si et seulement si
$P(\sX_\ell/\ell,\sF_\ell,Z_\ell)$ est vraie.

Une telle propri{\'e}t{\'e} $P$ {\'e}tant fix{\'e}e, revenons {\`a}
un morphisme de champs alg{\'e}briques $f:\sX\to\sS$ et un
point $s=[s_K]$. Pour tout repr{\'e}sentant
$s_K:\Spec(K)\to\sS$ de $s$, o{\`u} $K$ est un corps, on
note $\sX_K=\sX\times_\sS \Spec(K)$, $\sF_K$ la pr{\'e}image
de $\sF$ par la projection $\sX_K\to \sX$, $Z_K$ la
pr{\'e}image de $Z$ par l'application continue
$|\sX_K|\to |\sX|$. Compte tenu de l'hypoth{\`e}se sur $P$,
le fait que $P(\sX_K/K,\sF_K,Z_K)$ ait lieu est
ind{\'e}pendant du repr{\'e}sentant $s_K:\Spec(K)\to\sS$
choisi pour $s$. On dira alors que {\em la propri{\'e}t{\'e}
$P(\sX_s,\sF_s,Z_s)$ est vraie}. On notera qu'il s'agit
bien s{\^u}r d'un abus de langage, puisque nous n'avons
d{\'e}fini ni $\sX_s$, ni $\sF_s$, ni $Z_s$. {\`A} chaque fois
que nous utiliserons cette notation $P(\sX_s,\sF_s,Z_s)$,
il sera sous-entendu que la propri{\'e}t{\'e} $P$ est
ind{\'e}pendante du corps de base au sens ci-dessus (et
cette ind{\'e}pendance sera {\'e}vidente ou bien connue).
\end{appadoc}

\begin{appadoc}{Propri{\'e}t{\'e}s constructibles}
Dans les {\'e}nonc{\'e}s ci-dessous, la v{\'e}rification du
fait qu'une propri{\'e}t{\'e} $P$ ait lieu se ram{\`e}ne
toujours au cas o{\`u} $\sS$ est un sch{\'e}ma, apr{\`e}s
changement de base par une pr{\'e}sentation lisse $S\to\sS$.

Concernant la constructibilit{\'e}, on notera que si
$\sX$ est un champ alg{\'e}brique et $Z$ est une
partie de $|\sX|$, la propri{\'e}t{\'e} pour $Z$ d'{\^e}tre
constructible est locale sur $\sX$ pour la topologie lisse.
On a m{\^e}me une propri{\'e}t{\'e} beaucoup plus forte, puisque
$Z$ est constructible si et seulement si $u^{-1}(Z)$ est
constructible, pour n'importe quel morphisme surjectif
et ouvert $u:\sU\to\sX$.


\begin{quot}{{\bf (Chevalley)}}
Soit $f:\sX\to\sS$ un morphisme de pr{\'e}sentation finie de champs
alg{\'e}briques et soit $Z$ une partie constructible de $|\sX|$. Alors
$f(Z)$ est une partie constructible de $|\sS|$.
\end{quot}

Voir \cite{LMB}, th{\'e}or{\`e}me 5.9.4.

\begin{prop} \label{prop:EGA_95}
Soient $\sX\to\sS$ un morphisme de pr{\'e}sentation finie de champs
alg{\'e}briques et $Z,Z'$ deux parties localement constructibles de
$|\sX|$. Alors, les ensembles suivants sont localement constructibles
dans $|\sS|$~:
\begin{trivlist}
\itemn{i} l'ensemble des $s\in |\sS|$ tels que $Z_s\ne\emptyset$,
\itemn{ii} l'ensemble des $s\in |\sS|$ tels que $Z_s\subset Z'_s$
(resp. $Z_s=Z'_s$),
\itemn{iii} si $Z\subset Z'$, l'ensemble des $s\in |\sS|$ tels
que $Z_s$ est dense dans $Z'_s$,
\itemn{iv} l'ensemble des $s\in |\sS|$ tels que $Z_s$ est
ouvert (resp. ferm{\'e}, resp. localement ferm{\'e}) dans $|\sX_s|$.
\end{trivlist}
\end{prop}

Voir \cite{EGA}~IV.9.5.1, 9.5.2, 9.5.3, 9.5.4. Les preuves
de 9.5.1, 9.5.2 et de 9.5.4 (une fois
d{\'e}montr{\'e} 9.5.3) s'adaptent imm{\'e}diatement au cas des
champs alg{\'e}briques. La preuve de 9.5.3 pour un morphisme
de sch{\'e}mas $f:X\to S$ se ram{\`e}ne au cas o{\`u} $X$
est int{\`e}gre. Ensuite, on utilise un recouvrement de $X$
par des ouverts affines int{\`e}gres.
Dans le cas des champs alg{\'e}briques, l'utilisation d'un
recouvrement lisse de $\sX$ par des ouverts affines
int{\`e}gres fait tout aussi bien l'affaire.

\begin{prop}
Soient $\sX\to\sS$, $\sY\to\sS$ deux morphismes de
pr{\'e}sentation finie de champs alg{\'e}briques et $f:\sX\to\sY$
un $\sS$-morphisme. Alors, les ensembles suivants sont
localement constructibles dans $|\sS|$~:
\begin{trivlist}
\itemn{i} l'ensemble des $s\in |\sS|$ tels que $f_s$ est
une immersion,
\itemn{ii} l'ensemble des $s\in |\sS|$ tels que $f_s$ est
une immersion ferm{\'e}e,
\itemn{iii} l'ensemble des $s\in |\sS|$ tels que $f_s$ est
une immersion ouverte.
\end{trivlist}
\end{prop}

Voir \cite{EGA}~IV.9.6.1, (viii), (ix), (x).

La proposition suivante fait intervenir la notion de
{\em point associ{\'e}} d'un faisceau quasi-coh{\'e}rent sur
un champ alg{\'e}brique. La d{\'e}finition correcte est la
suivante~; elle est tir{\'e}e du paragraphe 2.2.6.3 de
Lieblich~\cite{Lie}, auquel nous renvoyons le lecteur
pour plus de d{\'e}tails.

\begin{defi} \label{defi:points_associes}
Soit $\sX$ un champ alg{\'e}brique localement noeth{\'e}rien et
$\sF$ un $\cO_\sX$-module quasi-coh{\'e}rent. On dit qu'un point
$x\in |\sX|$ est un {\em point associ{\'e}} de $\sF$ s'il
existe un sous-faisceau $\sG\subset \sF$ quasi-coh{\'e}rent
tel que $x\in\Supp(\sG)\subset\bar{\{x\}}$. On note
$\Ass(\sF)$ l'ensemble des points associ{\'e}s de $\sF$.
\end{defi}

\begin{prop} \label{prop:semi_continuite_de_sci}
Soient $\sS$ un champ alg{\'e}brique noeth{\'e}rien int{\`e}gre de point
g{\'e}n{\'e}rique $\eta$, $f:\sX\to\sS$ un morphisme de type fini, $\sF$
un $\cO_\sX$-module coh{\'e}rent. Si $\sF_\eta$ est sans cycle premier
associ{\'e} immerg{\'e}, il existe un voisinage $\sU$ de $\eta$ dans $\sS$
tel que pour tout $s\in |\sU|$, $\sF_s$ soit sans cycle premier
associ{\'e} immerg{\'e}.
\end{prop}

Voir \cite{EGA}~IV.9.7.6.
Le cas des champs alg{\'e}briques en d{\'e}coule en prenant une
pr{\'e}sentation lisse $\pi:X\to\sX$, car $\sF$ est sans cycle
premier associ{\'e} immerg{\'e} si et seulement si $\pi^*\sF$ est sans
cycle premier associ{\'e} immerg{\'e}.

\begin{theo} \label{theo:EGA_977}
Soit $\sX\to\sS$ un morphisme de pr{\'e}sentation finie de champs
alg{\'e}briques. Alors, les ensembles suivants sont
localement constructibles dans $|\sS|$~:
\begin{trivlist}
\itemn{i} l'ensemble des $s\in |\sS|$ tels que
$\sX_s$ est g{\'e}om{\'e}triquement irr{\'e}ductible,
\itemn{ii} l'ensemble des $s\in |\sS|$ tels que
$\sX_s$ est g{\'e}om{\'e}triquement connexe,
\itemn{iii} l'ensemble des $s\in |\sS|$ tels que
$\sX_s$ est g{\'e}om{\'e}triquement r{\'e}duit,
\itemn{iv} l'ensemble des $s\in |\sS|$ tels que
$\sX_s$ est g{\'e}om{\'e}triquement int{\`e}gre.
\end{trivlist}
\end{theo}

Voir \cite{EGA}~IV.9.7.7.
Pour adapter la d{\'e}monstration au cas des champs alg{\'e}briques,
on se ram{\`e}ne imm{\'e}diatement au cas o{\`u} $\sS$ est un sch{\'e}ma $S$
en prenant les images inverses par une pr{\'e}sentation lisse
$S\to\sS$.

Le principal point d{\'e}licat se situe alors dans le num{\'e}ro
1\textordmasculine\ de {\em loc. cit.} o{\`u} l'on construit un ouvert $W$
commun {\`a} $X$ et {\`a} un sch{\'e}ma de la forme
$Y=\Spec(A[T_1,\dots,T_{n+1}]/(F))$ (notations de {\em
  loc. cit.}). Lorsque $\sX$ est un champ alg{\'e}brique, nous proc{\'e}derons
comme suit. On consid{\`e}re une pr{\'e}sentation lisse quasi-compacte
$X_{0,\eta}\to\sX_\eta$ de la fibre de $\sX$ au-dessus du point
g{\'e}n{\'e}rique de $S$, et on pose $X_{1,\eta}=X_{0,\eta}\times_\sX
X_{0,\eta}$. On note $\delta_{0,i}$ ($1\le i\le n_0$)
(resp. $\delta_{1,j}$ $1\le j\le n_1$) les points g{\'e}n{\'e}riques de
$X_{0,\eta}$ (resp. de $X_{1,\eta}$) et $L_{0,i}$ (resp. $L_{1,j}$)
leurs anneaux locaux, qui sont des corps. On <<~d{\'e}ploie~>> chacun des
corps de fonctions $L_{a,i}$ en un sch{\'e}ma affine de la forme
$Y_{a,i}=\Spec(A[T_{a,i,1},\dots,T_{a,i,n+1}]/(F_{a,i}))$ par le
proc{\'e}d{\'e} de \cite{EGA}. Notons $Y_0$ resp. $Y_1$ le sch{\'e}ma
somme disjointe des $Y_{0,i}$ resp. des $Y_{1,i}$. Par construction
$Y_{0,\eta}$ resp. $Y_{1,\eta}$ est somme de sch{\'e}mas int{\`e}gres
de corps de fonctions rationnelles les $L_{a,i}$, et on dispose donc
de deux fl{\`e}ches $Y_{1,\eta}\dasharrow Y_{0,\eta}$ d{\'e}finies en
codimension $0$ (ce sont les restrictions des deux fl{\`e}ches
$X_{1,\eta}\to X_{0,\eta}$ aux points g{\'e}n{\'e}riques). Comme
$Y_0$, $Y_1$ sont de pr{\'e}sentation finie, quitte {\`a} les remplacer par
des ouverts, on peut supposer que ces deux fl{\`e}ches s'{\'e}tendent en des
morphismes $f,g:Y_1\to Y_0$. Comme au-dessus de~$\eta$ ces morphismes
co{\"\i}ncident g{\'e}n{\'e}riquement avec les deux fl{\`e}ches $X_{1,\eta}\to
X_{0,\eta}$ qui sont lisses, quitte {\`a} restreindre encore $Y_0$ et
$Y_1$ on peut supposer que $f,g$ sont lisses. Elles d{\'e}finissent donc
un groupo{\"\i}de lisse dont on note $\sY$ le champ alg{\'e}brique quotient. Ce
champ alg{\'e}brique joue le r{\^o}le tenu par $Y$ dans la preuve de
\cite{EGA}~IV.9.7.7.

On doit ensuite justifier qu'il existe un voisinage ouvert $U$ de
$\eta$ dans $S$ tel que $\sY_s$ reste g{\'e}om{\'e}triquement int{\`e}gre
pour tout $s\in U$. Or $\sY_s$ est g{\'e}om{\'e}triquement int{\`e}gre si et
seulement si $Y_{0,s}$ est g{\'e}om{\'e}triquement ponctuellement int{\`e}gre et
$Y_{1,s}\to Y_{0,s}\times Y_{0,s}$ est dominant. On obtient donc l'existence
d'un tel $U$ en utilisant \cite{EGA}~IV.9.7.4, 9.7.5 et 9.6.1(ii).

\begin{prop} \label{prop:EGA_979}
Soit $\sX\to\sS$ un morphisme de pr{\'e}sentation finie de champs
alg{\'e}briques. Alors les fonctions <<~nombre g{\'e}om{\'e}trique de
composantes connexes de $\sX_s$~>> et <<~nombre g{\'e}o\-m{\'e}trique de
composantes irr{\'e}ductibles de $\sX_s$~>> sont localement constructibles.
\end{prop}

Voir \cite{EGA}~IV.9.7.9. La preuve utilise IV.9.7.8 et IV.9.7.1,
dont les {\'e}nonc{\'e}s et les preuves s'adaptent sans
modification pour les champs alg{\'e}briques. Pour ne pas alourdir
inutilement le texte, nous ne les recopions pas ici.
\end{appadoc}

\section{Puret{\'e} pour les champs alg{\'e}briques}
\label{annexe:purete}

Si $X\to S$ est un morphisme de sch{\'e}mas localement de type
fini et $\sM$ est un $\cO_X$-module quasi-coh{\'e}rent de
pr{\'e}sentation finie, la notion de {\em puret{\'e} de $\sM$
relativement {\`a} $S$} est d{\'e}finie dans \cite{RG},~3.3.3.
Si $\sM$ est plat sur $S$, cette d{\'e}finition est locale
pour la topologie plate sur $S$ (\cite{RG}, 3.3.7). En revanche,
elle n'est
pas locale sur $X$, m{\^e}me pour la topologie de Zariski.
En cons{\'e}quence, l'extension de cette notion au cas o{\`u}
$\sX\to\sS$ est un morphisme de champs alg{\'e}briques ne peut
se faire simplement en prenant une pr{\'e}sentation lisse de
$\sX$. On en revient donc {\`a} la d{\'e}finition originale,
passant par les points associ{\'e}s
(d{\'e}finition~\ref{defi:points_associes}).

\begin{appdefi}
Soit $f:\sX\to\sS$ un morphisme localement de type fini de
champs alg{\'e}briques et soit $\sM$ un $\cO_\sX$-module
quasi-coh{\'e}rent de pr{\'e}sentation finie, plat sur $\sS$.
\begin{trivlist}
\itemn{i} Supposons que $\sS$ est un sch{\'e}ma local
hens{\'e}lien $S$ de point ferm{\'e} $s_0$. On dit que $\sM$
est pur relativement {\`a} $S$ si pour tout $x\in |\sX|$,
qui est un point associ{\'e} du $\cO_{\sX_s}$-module $\sM_s$,
o{\`u} $s=f(x)$, l'adh{\'e}rence de $x$ dans $|\sX|$
rencontre $\sX_{s_0}$.
\itemn{ii} Supposons que $\sS$ est repr{\'e}sentable par
un espace alg{\'e}brique $S$. On dit que $\sM$ est pur
relativement {\`a} $S$ si pour tout $s\in S$, de
hens{\'e}lis{\'e} $(S^h,s^h)$, le module $\sM\times_S S^h$
est pur relativement {\`a} $S^h$.
\itemn{iii} En g{\'e}n{\'e}ral, on dit que $\sM$ est pur
relativement {\`a} $\sS$ si $\pi^*\sM$ est pur
relativement {\`a} $S$, pour une (et donc toute)
pr{\'e}sentation lisse $\pi:S\to\sS$.
\end{trivlist}
\end{appdefi}

\begin{applemm} \label{lemm:connexite_et_fibre_generique}
Soit $R=(R,K,k,\pi)$ un anneau de valuation discr{\`e}te
et $\sX$ un $R$-champ alg{\'e}brique localement
de type fini, plat, {\`a} fibre sp{\'e}ciale r{\'e}duite.
Alors, $\sX$ est connexe si et seulement si $\sX_K$ est
connexe, et $\sX$ est int{\`e}gre si et seulement si
$\sX_K$ est int{\`e}gre.
\end{applemm}

\begin{demo}
Posons $B=H^0(\sX,\cO_\sX)$, on a $B_K=H^0(\sX_K,\cO_{\sX_K})$.
Pour la premi{\`e}re assertion, il suffit de montrer que les
idempotents de $B$ et $B_K$ sont les m{\^e}mes, et donc de montrer
que les idempotents de $B_K$ sont dans $B$. Soit $e\in B_K$ tel
que $e^2=e$. Si $e\not\in B$, il existe une {\'e}criture
$e=\pi^{-n}f$ avec $n\ge 1$ et $f\in B\setminus\pi B$. Comme
$f^2=(\pi^ne)^2=\pi^n f$ et que $B_k$, en tant que sous-anneau
de $H^0(\sX_k,\cO_{\sX_k})$, est r{\'e}duit par hypoth{\`e}se,
on trouve que $f$ est nul modulo $\pi$, contradiction.
Pour la seconde assertion, il ne reste qu'{\`a} montrer que $\sX$
est localement int{\`e}gre si et seulement si $\sX_K$ l'est. Or par
platitude cela est clair si $\sX$ est un sch{\'e}ma, et on se ram{\`e}ne
{\`a} ce cas en utilisant une pr{\'e}sentation lisse $X\to \sX$.
\end{demo}

\begin{applemm} \label{lemm:revetement_par_U_pur}
Soit $R$ un anneau de valuation discr{\`e}te hens{\'e}lien et
$\sX$ un $R$-champ alg{\'e}brique de type fini, plat et pur.
Alors, il existe un $R$-sch{\'e}ma affine $U$ avec une
$R$-alg{\`e}bre de fonctions de type fini, libre comme $R$-module,
et un morphisme lisse $R$-universellement sch{\'e}matiquement
dominant $U\to\sX$. En particulier $H^0(\sX,\cO_\sX)$ est un
$R$-module libre.
Si $\sX$ est {\`a} fibre sp{\'e}ciale irr{\'e}ductible, on peut
supposer que $U$ est r\'eunion disjointe d'un nombre fini d'ouverts
affines {\`a} fibre sp{\'e}ciale irr{\'e}ductible et alg{\`e}bre
de fonctions de type fini libre comme $R$-module.
\end{applemm}

\begin{demo}
Fixons une pr{\'e}sentation lisse $X\to\sX$. En chaque point
$x\in X$ de la fibre sp{\'e}ciale, choisissons un voisinage ouvert
affine $U_x$. Si $\sX_k$ est irr{\'e}ductible, on peut choisir
$U_x$ {\`a} fibre sp{\'e}ciale irr{\'e}ductible. Quitte {\`a}
r{\'e}tr{\'e}cir $U_x$, on peut supposer de plus que son anneau
de fonctions
est s{\'e}par{\'e} pour la topologie $\pi$-adique (voir
\cite{Ro2}, lemma~2.1.11). D'apr{\`e}s Raynaud et Gruson, un
tel anneau est libre comme $R$-module (voir \cite{Ro2},
lemma~2.1.7). Comme $\sX$ (et
donc aussi $\sX_k$) est quasi-compact, un nombre fini des ouverts
$U_x$ recouvre $\sX_k$. Soit $U$ la somme disjointe de ces ouverts.
Comme $\sX$ est pur, aucun cycle premier associ{\'e} de $\sX$
n'est inclus dans $\sX_K$. Un tel cycle premier associ{\'e} est alors
inclus dans l'image de $U$, de sorte que le morphisme $U\to\sX$ est
sch{\'e}matiquement dominant. Comme $U_k\to\sX_k$ l'est aussi,
il s'ensuit que $U\to \sX$ est universellement sch{\'e}matiquement
dominant (par l'argument de \cite{Ro2},~lemma~2.1.9).
Il s'ensuit qu'on a
une injection $H^0(\sX,\cO_\sX)\into H^0(U,\cO_U)$. Comme un
sous-module libre d'un module libre sur un anneau de valuation
discr{\`e}te est libre, ceci prouve le lemme.
\end{demo}

Le th{\'e}or{\`e}me 2.2.1 de \cite{Ro2} se g{\'e}n{\'e}ralise
imm{\'e}diatement~:

\begin{apptheo} \label{theo:CC_CI_pour_les_champs}
Soit $f:\sX\to\sS$ un morphisme de champs alg{\'e}briques
de pr{\'e}sentation finie, plat et pur, et soit $n\ge 1$ un
entier. Alors, les ensembles suivants sont ouverts dans
$|\sS|$~:
\begin{trivlist}
\itemn{i} l'ensemble des $s\in |\sS|$ tels que $\sX_s$
est g{\'e}om{\'e}triquement r{\'e}duit,
\itemn{ii} l'ensemble des $s\in |\sS|$ tels que la
fibre g{\'e}om{\'e}trique $\sX_{\bar s}$ est r{\'e}duite avec
au plus $n$ composantes connexes,
\itemn{iii} l'ensemble des $s\in |\sS|$ tels que la
fibre g{\'e}om{\'e}trique $\sX_{\bar s}$ est r{\'e}duite avec
au plus $n$ composantes irr{\'e}ductibles.
\end{trivlist}
\end{apptheo}

\begin{demo}
On adapte la preuve de~\cite{Ro2},~th.2.2.1. En faisant le
changement de base par une pr{\'e}sentation lisse $S\to \sS$,
on se ram{\`e}ne au cas o{\`u} $\sS$ est un sch{\'e}ma $S$. On se
ram{\`e}ne ensuite au cas o{\`u} $S$ est affine, noeth{\'e}rien.
D'apr{\`e}s~\ref{theo:EGA_977}~(iii) et \ref{prop:EGA_979},
les ensembles qui nous int{\'e}ressent sont constructibles
dans~$S$. On se ram{\`e}ne alors au cas o{\`u} $S$ est le spectre
d'un anneau de valuation discr{\`e}te $R=(R,K,k,\pi)$, complet
{\`a} corps r{\'e}siduel alg{\'e}briquement clos, puis, en prenant
la cl{\^o}ture int{\'e}grale de $R$ dans une extension finie de
$K$, on se ram{\`e}ne {\`a} prouver que $\sX_K$ est r{\'e}duit
(resp. poss{\`e}de au plus $n$ composantes connexes, resp.
poss{\`e}de au plus $n$ composantes irr{\'e}ducitbles) d{\`e}s que
$\sX_k$ poss{\`e}de la m{\^e}me propri{\'e}t{\'e}.

\medskip

\no (i) On prouve que $\sX_K$ est r{\'e}duit exactement comme dans
\cite{Ro2}, th.~2.2.1, en rempla{\c c}ant l'ouvert universellement
sch{\'e}matiquement dominant $U$ r{\'e}union d'ouverts affines purs
utilis{\'e} dans {\em loc. cit.} par le morphisme $U\to \sX$ fourni
par le lemme~\ref{lemm:revetement_par_U_pur}.

\medskip

\no (ii) La preuve de \cite{Ro2}, th.~2.2.1 est valable sans
modification.

\medskip

\no (iii) Le d{\'e}but de la preuve de \cite{Ro2}, th.~2.2.1 est
valable sans modification, jusqu'au moment o{\`u} l'on fait appel
{\`a} un recouvrement ouvert sch{\'e}matiquement dominant par des
affines purs. On remplace le recours {\`a} ce recouvrement par
l'utilisation du morphisme $U\to \sX$ fourni par le
lemme~\ref{lemm:revetement_par_U_pur}, o{\`u} $U$ est somme
disjointe d'ouverts affines $U_x$ {\`a} fibre sp{\'e}ciale int{\`e}gre
et d'anneau de fonction libre comme $R$-module.
\end{demo}


\begin{thebibliography}{99}
\bibitem[AOV]{AOV} {\sc D. Abramovich, M. Olsson, A. Vistoli},
{\it Tame stacks in positive characteristic}, Ann. Inst. Fourier
(Grenoble) 58 (2008), no. 4, 1057--1091. 
\bibitem[Ar1]{Ar1} {\sc M. Artin}, {\it Versal deformations and algebraic
stacks}, Invent. Math. 27 (1974), 165--189.
\bibitem[BR]{BR} {\sc J. Bertin, M. Romagny}, {\em Champs de Hurwitz},
pr{\'e}publication disponible {\`a} l'adresse
http:/$\!$/people.math.jussieu.fr/\raisebox{-1.5mm}{\textasciitilde}romagny/.
\bibitem[Bro]{Bro} {\sc S. Brochard}, {\it Foncteur de Picard d'un champ
alg{\'e}brique}, Math. Ann. 343 (2009), 541--602.
\bibitem[DM]{DM} {\sc P. Deligne, D. Mumford},
{\it The irreducibility of the space of curves of given genus},
Publ. Math. IH{\'E}S~36 (1969), 75--109. 
\bibitem[EGA]{EGA} {\sc J.~Dieudonn{\'e}, A.~Grothendieck},
{\it {\'E}l{\'e}ments de G{\'e}om{\'e}trie Alg{\'e}brique II, III, IV},
Publ. Math. IH{\'E}S~8 (1961), 17 (1963), ~24 (1965), 28 (1966),
32 (1967).
\bibitem[Ei]{Ei} {\sc D. Eisenbud}, {\it Commutative algebra with a view
  toward algebraic geometry}, Graduate Texts in Math., Springer-Verlag (1995).
\bibitem[Kn]{Kn} {\sc D. Knutson}, {\it Algebraic spaces},
Lecture Notes in Mathematics, Vol.~203, Springer-Verlag, 1971.
\bibitem[La]{La} {\sc Y. Laszlo}, {\it Linearization of group stack
actions and the Picard group of the moduli of $\SL_r/\mu_s$-bundles
on a curve}, Bull. Soc. Math. France 125 (1997), no. 4, 529--545.
\bibitem[Lie]{Lie} {\sc M. Lieblich}, {\it Moduli of twisted
sheaves}, Duke Math. J. 138 (2007), no.~1, 23--118.
\bibitem[LMB]{LMB} {\sc G. Laumon, L. Moret-Bailly}, {\it Champs
alg{\'e}briques}, Ergebnisse der Mathematik und ihrer Grenzgebiete (3.~Folge)
no~39, Springer-Verlag, Berlin, 2000.
\bibitem[MSSV]{MSSV} {\sc K. Magaard, T. Shaska, S. Shpectorov,
H. V\"{o}lklein}, {\it The locus of curves with prescribed
automorphism group}, Communications in arithmetic fundamental
groups (Kyoto, 1999/2001).
S\={u}rikaisekikenky\={u}sho K\={o}ky\={u}roku no.~1267 (2002), 112--141.
\bibitem[Ol]{Ol} {\sc M.~Olsson}, {\it On proper coverings of Artin
stacks}, Adv. Math. 198 (2005), no.~1, 93--106. 
\bibitem[RG]{RG} {\sc M. Raynaud, L. Gruson}, {\it Crit{\`e}res de
platitude et de projectivit{\'e}. Techniques de <<~platification~>>
d'un module}, Invent. Math. 13 (1971), 1--89.
\bibitem[Ro1]{Ro1} {\sc M. Romagny}, {\it Group actions on stacks and
applications}, Michigan Math. J. 53 (2005), no.~1, 209--236.
\bibitem[Ro2]{Ro2} {\sc M. Romagny}, {\it Effective models of group
schemes}, http://arxiv.org/abs/0904.3167.
\bibitem[TV]{TV} {\sc B. To\"{e}n, G. Vezzosi}, {\it Homotopical
algebraic geometry~II. Geometric stacks and applications},
Mem. Amer. Math. Soc. 193 (2008),  no.~902, x+224 pp.
\end{thebibliography}
\end{document}